\documentclass[smallextended]{svjour3}
\usepackage{amssymb}
\usepackage{amsmath}
\usepackage{dsfont}
\usepackage{epsfig}
\usepackage{multicol}
\usepackage{enumerate}
\usepackage{caption}
\usepackage{psfrag}
\usepackage{xcolor}
\usepackage{graphicx}
\usepackage{authblk}
\usepackage{fullpage}
\usepackage[authoryear,sectionbib,sort]{natbib}
\usepackage{float}
\usepackage[caption = false]{subfig}
\input epsf

\numberwithin{equation}{section}

\newtheorem{hyps}{Hypotheses}
\newtheorem{hyp}[hyps]{Hypothesis}

\newcommand{\be}{\begin{equation}}
\newcommand{\ee}{\end{equation}}
\newcommand{\baa}{\begin{array}}
\newcommand{\eaa}{\end{array}}
\newcommand{\ba}{\begin{eqnarray}}
\newcommand{\ea}{\end{eqnarray}}

\newcommand{\R}{{\mathbb R}}
\newcommand{\N}{{\mathbb N}}

\newcommand{\bee}{\begin{equation*}}
\newcommand{\eee}{\end{equation*}}
\newcommand{\bc}{\begin{cases}}
\newcommand{\ec}{\end{cases}}

\journalname{}

\begin{document}

\title{Climate change and integrodifference equations in a stochastic environment}
\author{Juliette Bouhours        \and
        Mark A. Lewis %etc.
}


\institute{J. Bouhours \at
              Department of Mathematical and Statistical Sciences, University of Alberta, Edmonton T6G2G1, Canada \\
              \email{bouhours@ualberta.ca}           %  \\
           \and
           M. A. Lewis \at
              Department of Mathematical and Statistical Sciences, University of Alberta, Edmonton T6G2G1, Canada
}

\date{Received: date / Accepted: date}

\maketitle

\begin{abstract}
Climate change impacts population distributions, forcing some species to migrate poleward if they are to survive and keep up with the suitable habitat that is shifting with the temperature isoclines.  Previous studies have analyzed whether populations have the capacity to keep up with shifting temperature isoclines, and have mathematically determined the combination of growth and dispersal that is needed to achieve this.  However, the rate of isocline movement can be highly variable, with much uncertainty associated with yearly shifts. The same is true for population growth rates.  Growth rates can be variable and uncertain, even within suitable habitats for growth. In this paper we reanalyze the question of population persistence in the context of the uncertainty and variability in isocline shifts and rates of growth. Specifically, we employ a stochastic integrodifference equation model on a patch of suitable habitat that shifts poleward at a random rate. We derive a metric describing the asymptotic growth rate of the linearized operator of the stochastic model.  This metric yields a threshold criterion for population persistence. We demonstrate that the variability in the yearly shift and in the growth rate has a significant negative effect on the persistence in the sense that it decreases the threshold criterion for population persistence. 
Mathematically, we show how the persistence metric can be connected to the principal eigenvalue problem for a related integral operator, at least for the case where isocline shifting speed is deterministic. Analysis of dynamics for the case where the dispersal kernel is Gaussian leads to the existence of a critical shifting speed, above which the population will go extinct, and below which the population will persist. This leads to clear bounds on rate of environmental change if the population is to persist.  Finally we illustrate our different results for butterfly population using numerical simulations and demonstrate how increased variances in isocline shifts and growth rates translate into decreased likelihoods of persistence. 
\end{abstract}

\maketitle
\section{Introduction}
The consequences of climate change on population abundance and distribution have been widely investigated for the last two decades. One of these consequences is a modification of range distributions. Indeed, we know that, for a large variety of vertebrate and invertebrate species, climate change induces range shift toward the poles or higher altitudes, contraction or expansion of the habitat and habitat loss (\citealt{PY03,Hetal06,Metal14,P06,Letal08}, amongst other). Mathematical models and simulations that include climate change have predicted an effect of climate change on range distribution through habitat migration, habitat reduction and expansion and habitat loss \citep{PY03,Ni00,Huetal15,MM00,Petal11,Petal12,Hetal13}. In this paper we are interested in understanding, with the aid of a mechanistic model, the effect of shifting range on population persistence when the yearly range shifts and the population growth rates are stochastic.
 
Mechanistic models have been used to study the effect of shifting boundaries on population persistence, by considering a suitable habitat, which is a bounded domain where the population can grow, that is shifted toward the pole at a forced speed $c>0$. \citet{PL} used a reaction-diffusion system with a moving suitable habitat to investigate the effect of climate change and shifting boundaries on population persistence when two populations compete with one another. \citet{BDNZ} used a similar equation, in a scalar framework to study the persistence property of one population facing shifting range, and characterised persistence as it depends on the shifting speed. More recent papers also investigate the effect of shifting range on population persistence of single populations \citep{Richteretal12,Lerouxetal13,Lietal14}. 
In all these reaction-diffusion equations are used to model the temporal evolution of the density $u$ of a population in space and time. That is, individuals are assumed to disperse and grow simultaneously. In these models dispersal is local in the sense that the population disperses to its closest neighbourhood, in a diffusive manner.
 
Another approach to modelling the temporal evolution of the density of a population, is to consider populations that disperse and reproduce successively, and to allow for nonlocal dispersal. In this case, integrodifference equations are the appropriate model for the dynamics of the density $u$. integrodifference equations, introduced by \citet{KS86} to model discrete-time growth-dispersal, assume time ($t=0, 1\ldots$) to be discrete, and space $\xi\in\Omega$ to be continuous. From one generation to the next the population grows, according to a nonlinear growth function $f(u)$ and then disperses according to a dispersal kernel $K$, so
\begin{equation}
\label{eq:IDE}
u_{t+1}(\xi)=\int_\Omega K(\xi,\eta)f(u_t(\eta))d\eta, \quad t\in\N.\\
\end{equation}
Strictly, the dispersal kernel $K(\xi,\eta)$ is a probability density function describing the chance of dispersal from $\eta$ to $\xi$.

\begin{figure}
\begin{center}
\subfloat[]{\includegraphics[scale=0.4]{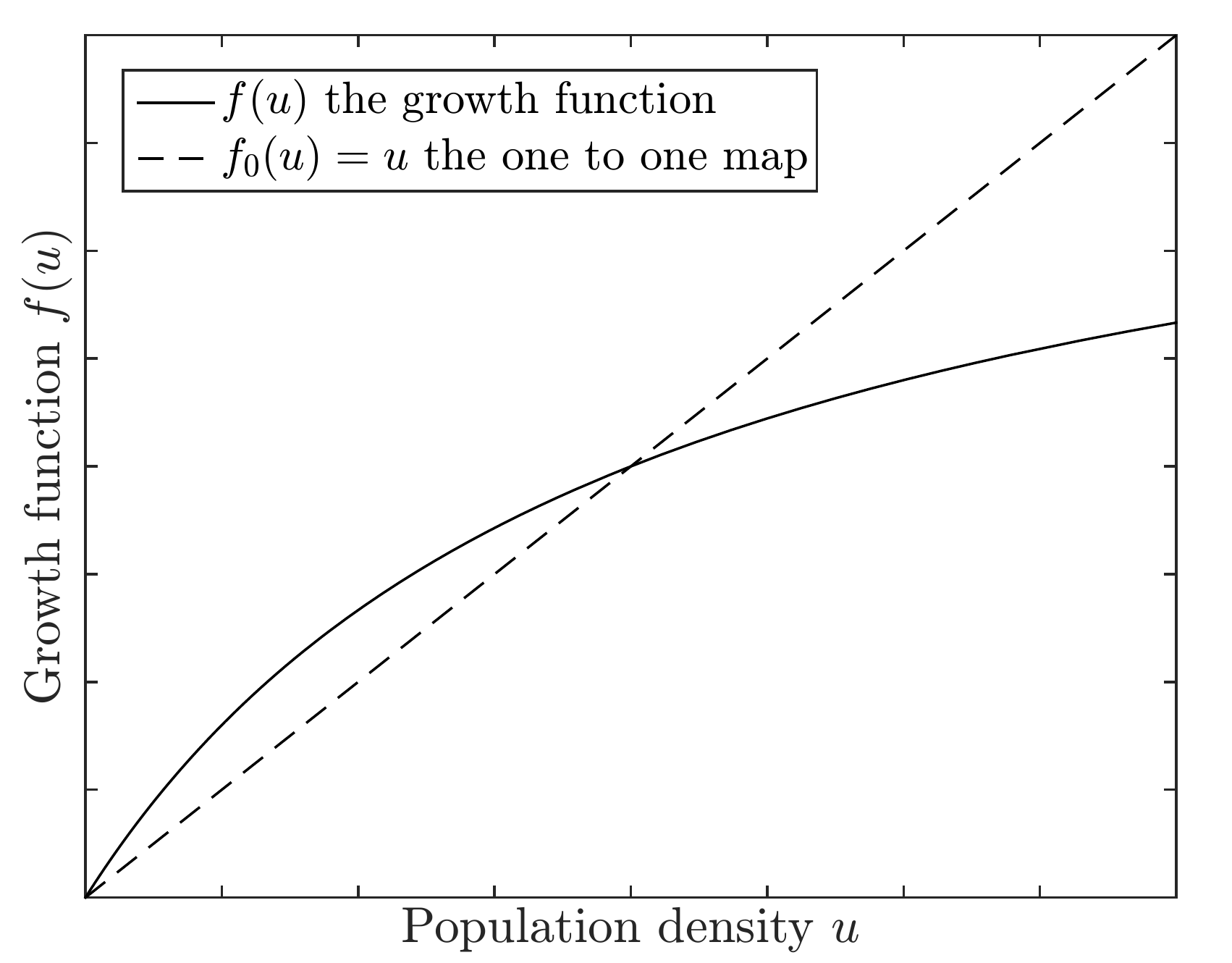}}\quad
\subfloat[]{\includegraphics[scale=0.4]{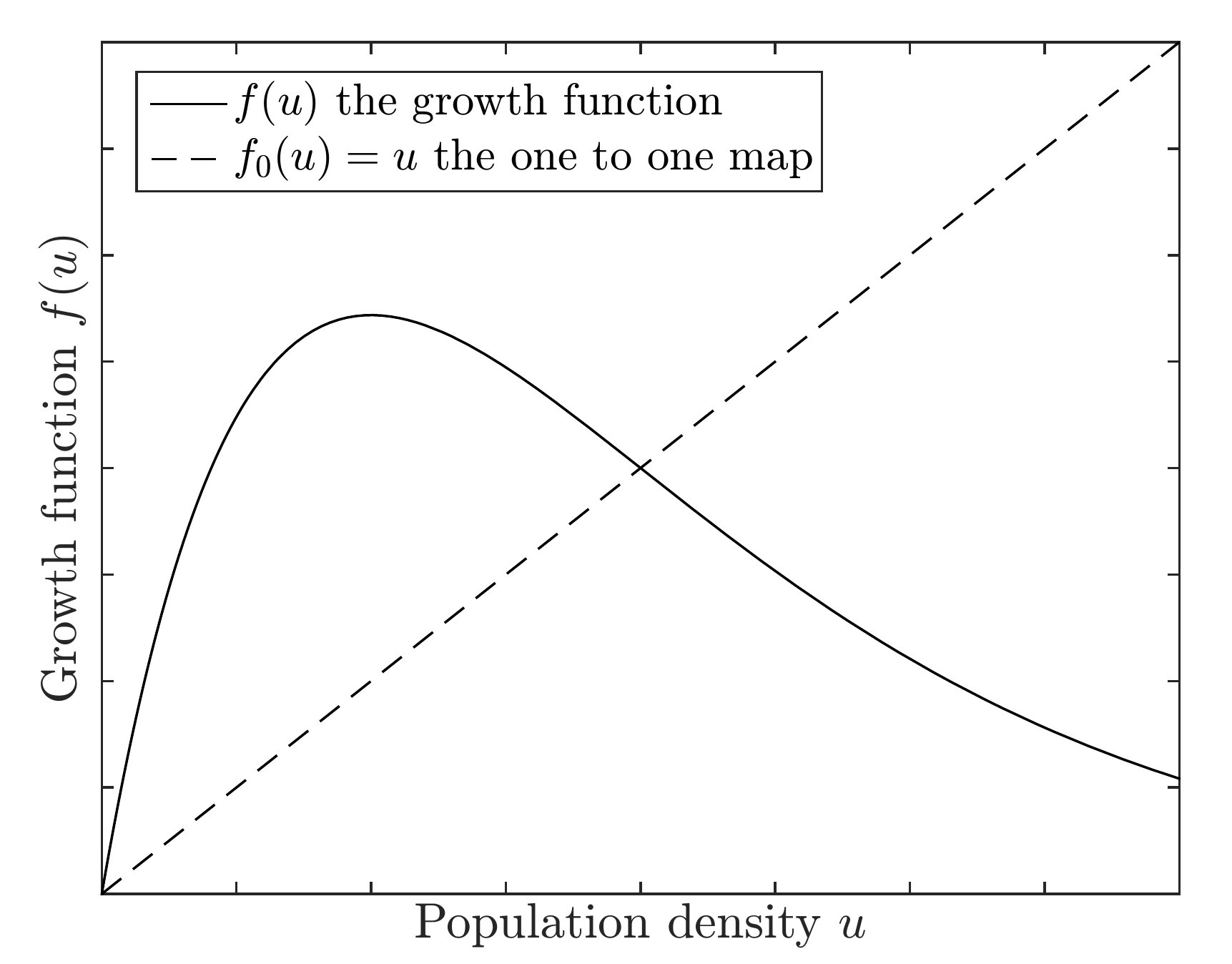}}
\caption{Compensatory (a) and over-compensatory (b) growth with one positive fixed point}\label{fgrowth}
\end{center}
\end{figure}

We consider a self-regulating population, with negative density-dependence so the slope of the growth is assumed to be monotically decreasing with $f(u)>u$ for $0<u<C$ and $f(u)<u$ for $u>C$, where $C>0$ is the carrying capacity. As we consider a population that is not subject to an Allee effect, the standard assumption on the growth function is that the geometric growth rate is the largest at lowest density, that is $f(u)/u$ achieves its supremum as $u$ approaches 0. We denote $r=\lim_{u\to 0^+} f(u)/u=f'(0)$.  When we wish to explicitly distinguish between populations with different geometric growth rates we modify our notation, replacing $f(u)$ by $f_r(u)$. Within this framework, we consider two types of growth dynamics: compensatory and over-compensatory. The compensatory growth dynamics are monotonic with respect to density $u$ whereas the overcompensatory growth dynamics have a characteristic ``hump'' shape.
(Figure \ref{fgrowth}). 

The eigenvalue problem associated with the linearisation of \eqref{eq:IDE} about $u\equiv 0$ is 
\be
\lambda \phi(\xi)=r\int_\Omega K(\xi,\eta)\phi(\eta)d\eta
\ee
and persistence of the population $u_t(\xi)$ depends upon whether $\lambda$ falls above or below one \citep{KS86, VL97,LL04}. An approximate method for calculating $\lambda$ employs the so-called dispersal success approximation \citep{VL97}.
Without loss of generality one can assume that $\int_{\Omega}\phi(\eta)d\eta=1$ and integrating the previous equation on $\Omega$ we get
\be\lambda=r\int_{\Omega}\int_{\Omega}K(\xi,\eta)\phi(\eta)d\eta=r\int_{\Omega}s(\eta)\phi(\eta)d\eta\ee
where $s(\eta)=\int_{\Omega}K(\xi,\eta)d\xi$ is the dispersal success. This function represents the probability for an individual located at $\eta$ to disperse to a point within the domain.
The so-called dispersal success approximation $\phi(\eta)\approx \frac{1}{|\Omega|}$ thereby allowing the principal eigenvalue $\lambda$ to be be estimated by
\be\label{eq:DispSuccAppr}\overline{\lambda}=\frac{r}{|\Omega|}\int_{\Omega}s(\eta)d\eta\ee
\citep{VL97}. This approximation gives $\overline{\lambda}$ as the growth rated times the estimated proportion of individuals that stay within the suitable habitat from one generation to the next. A modified dispersal success approximation recently introduced by \citet{RBM15} improves upon the dispersal success approximation assumption that the population is uniformly distributed within the favourable environment. \citet{RBM15} introduced a modified approximation that weights the dispersal success values by the proportion of the population at each point. They defined the modified dispersal success approximation by
\be\label{eq:ModDispSuccAppr}
\widehat{\lambda}= \frac{r}{|\Omega|}\int_{\Omega}\left(\frac{s(\eta)}{\overline{\lambda}}\right)s(\eta)d\eta\ee
and showed that this gave a better approximation to the eigenvalue.  We will employ both versions of the dispersal success approximation \eqref{eq:DispSuccAppr} and the modified dispersal success approximation \eqref{eq:ModDispSuccAppr} in our calculations later in this paper.

So far we have considered only the dependence of the growth on the density $u$. In the framework for climate change, the population can grow differently depending on where it is located with respect to space and time. To take this into account we introduce a suitability function, $g_t(\eta)$ ($0\leq g_t \leq 1$), which depends on space and time and multiplies the growth map $f$. As we consider populations whose suitable habitat shifts toward the pole, we choose a particular form for suitability function, $g_t(\eta)=g_0(\eta-s_t)$, where $g_0$ is the initial suitability function in the absence of climate change and $s_t$ is a parameter standing for the center of the suitable habitat.  The simple case where the habitat shifts at a constant speed $c$ is given by $s_t=ct$.  However, as we describe below, it is also possible to allow $s_t$ to vary randomly about $ct$.

\begin{figure}
\begin{center}
\includegraphics[scale=0.5]{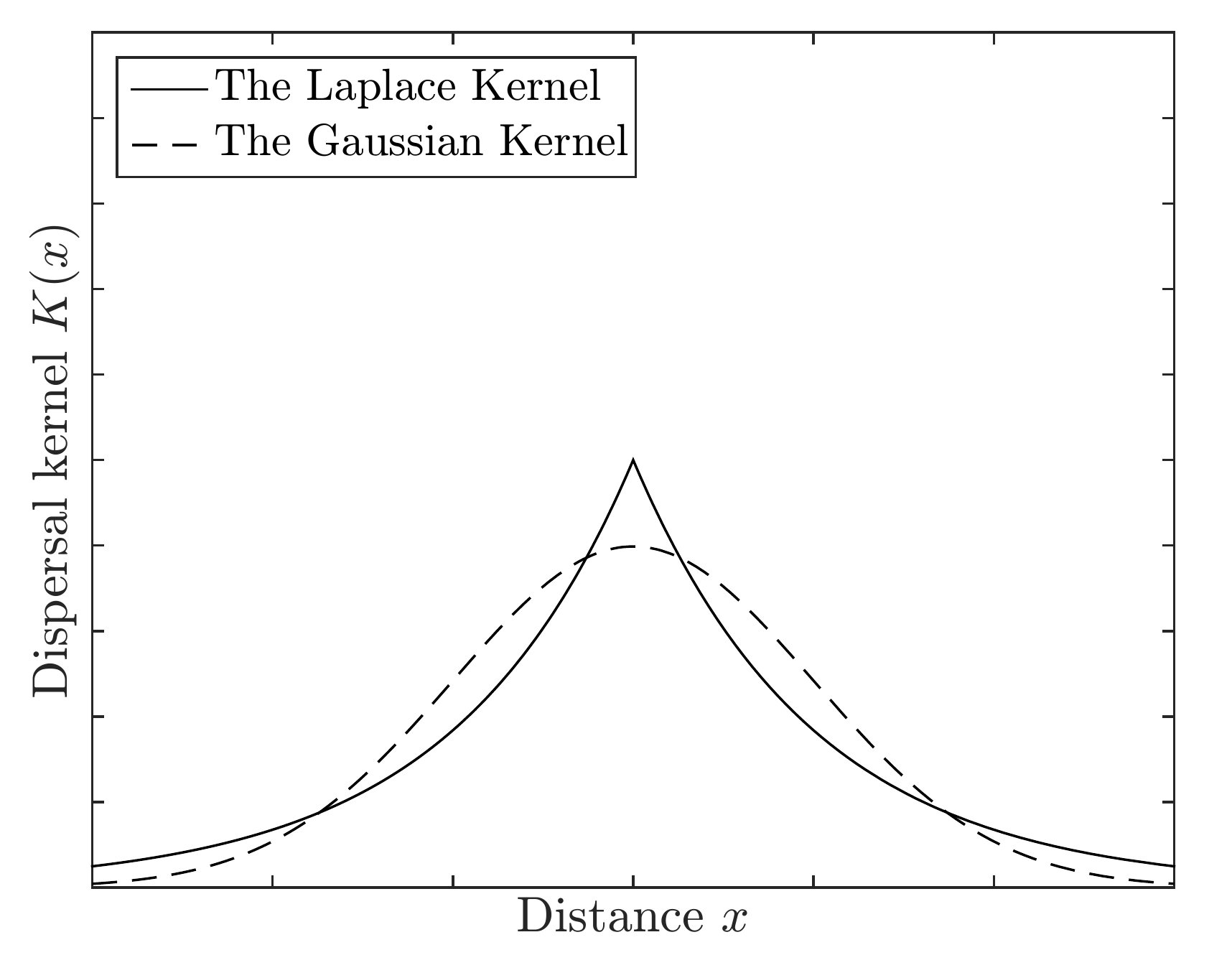}\quad
\caption{Laplace and Gaussian dispersion kernel $x\mapsto K(x)$ centered at 0}\label{fig:kernel}
\end{center}
\end{figure}

Following the growth stage, the population disperses in space according to the dispersal kernel $K$. The kernel $K$ is assumed to be positive everywhere: that is, the probability of dispersing to any point in space is always positive. Dispersal kernels are typically assumed to depend only on the signed distance between two points in space, i.e. only on the dispersal location relative to the source location.  If the population has no preferred direction of dispersal, the kernel is symmetric, depending only upon the distance between source and dispersal locations. This is not the case, however, in rivers where the population is subjected to a stream flow for example, or in environments with a prevailing wind direction that affects dispersal. In this paper we consider the Gaussian and the Laplace as examples of typical symmetric dispersal kernels (Figure \ref{fig:kernel})

In this framework of discrete-time growth-dispersal models, \citet{ZK1} investigated the effect of climate change and shifting range on the persistence of the population and highlighted the possible existence of a critical shifting speed for persistence. More recently several works investigated the effect of shifting range in more general discrete time growth-dispersal models \citep{ZK2, HZHK,PK15}. 

Recent research reports an increase in the environmental stochasticity in population dynamics, partly due to climate change and its effect on the frequency and the intensity of extreme climatic events covering large areas of the globe \citep{Saltzetal06,IPCC07,Kreylingetal11}. It is also known that the projected consequences for population ranges vary, depending on the different scenarios related to climate change (see, for example, \citet{IPCC14}). 

In this paper we focus on the effect of environmental stochasticity on population persistence in the presence of a shifting range, using a stochastic population model. The development of stochastic population models in population ecology was initially motivated by the study of the effect of environmental stochasticity on population dynamics and on the large time behaviour of the population \citep{May73, T77}.

\begin{figure}
\begin{center}
\includegraphics[scale=0.5]{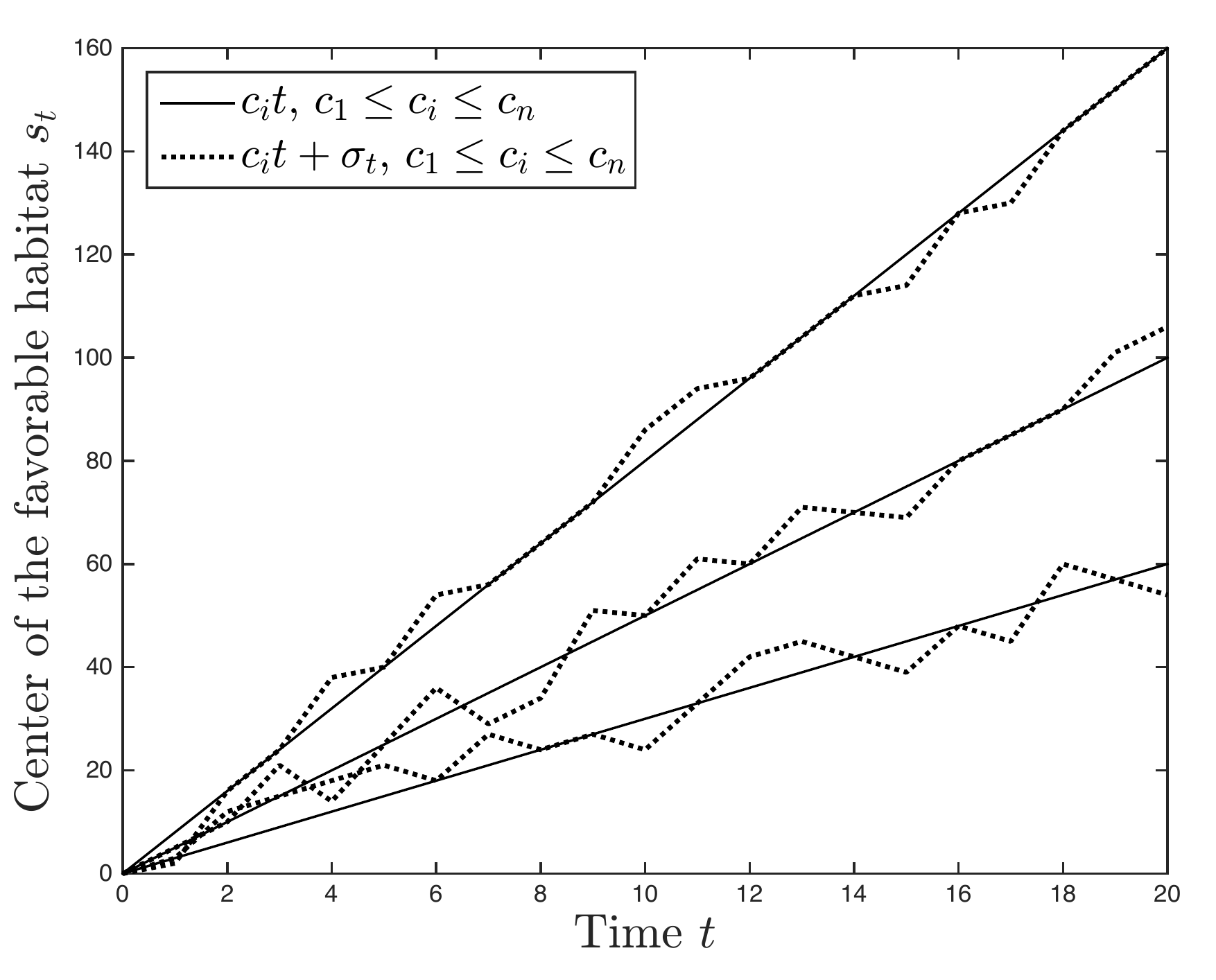}
\caption{Position of the center of the suitable habitat depending of the generation $t$, the chosen shifting speed $c_i\in\{c_1,\dots,c_n\}$ and the random process $\sigma_t$. In this figure there are three possible $c$ and for each $c_i\in\{c_1,c_2,c_3\}$, the center of the suitable habitat at generation $t$ is located at $c_i t+\sigma_t$, $(\sigma_t)_t$ independently identically distributed}\label{fig:CSigma}
\end{center}
\end{figure}

\begin{figure}
\begin{center}
\subfloat[]{\includegraphics[scale=0.4]{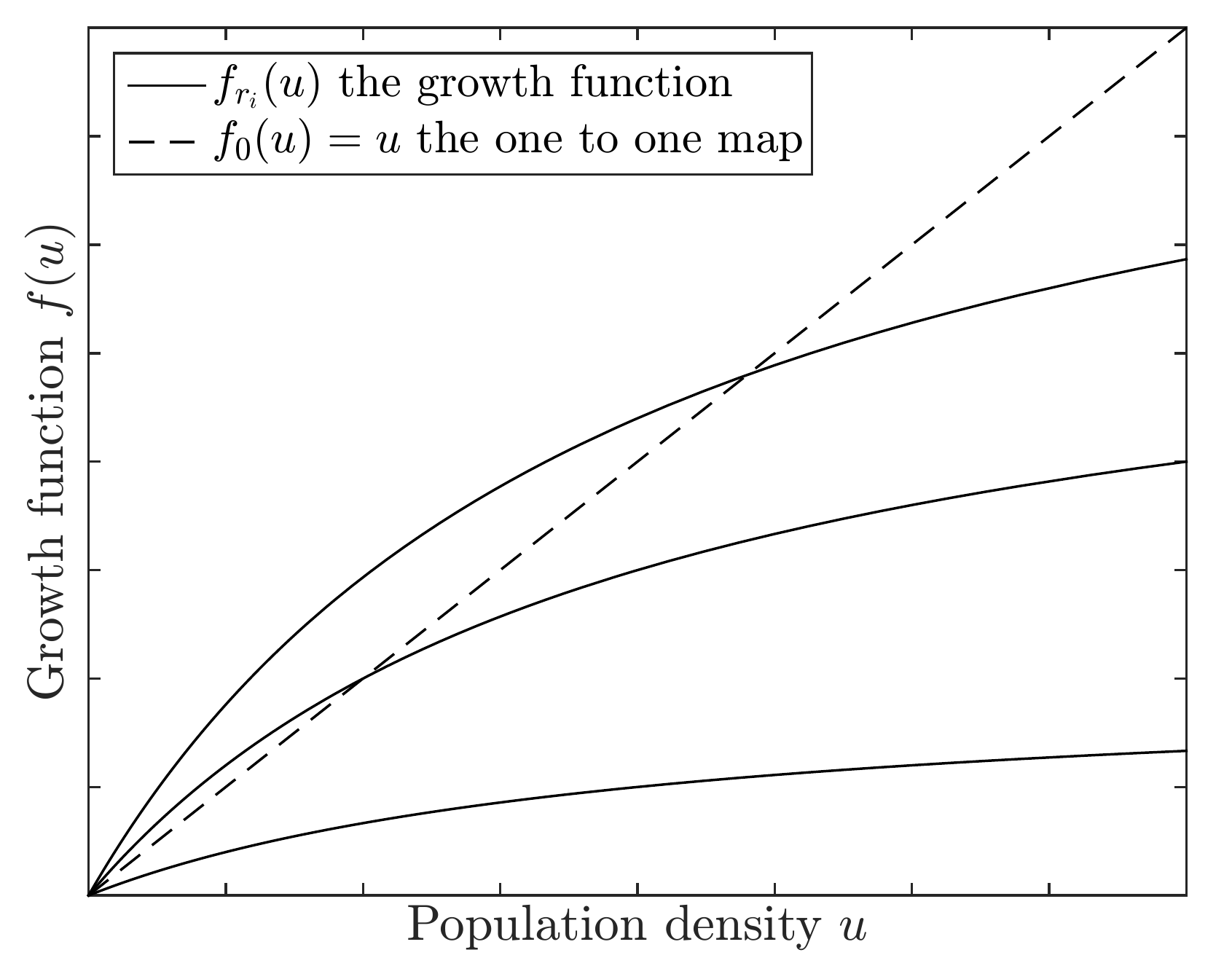}}\quad
\subfloat[]{\includegraphics[scale=0.4]{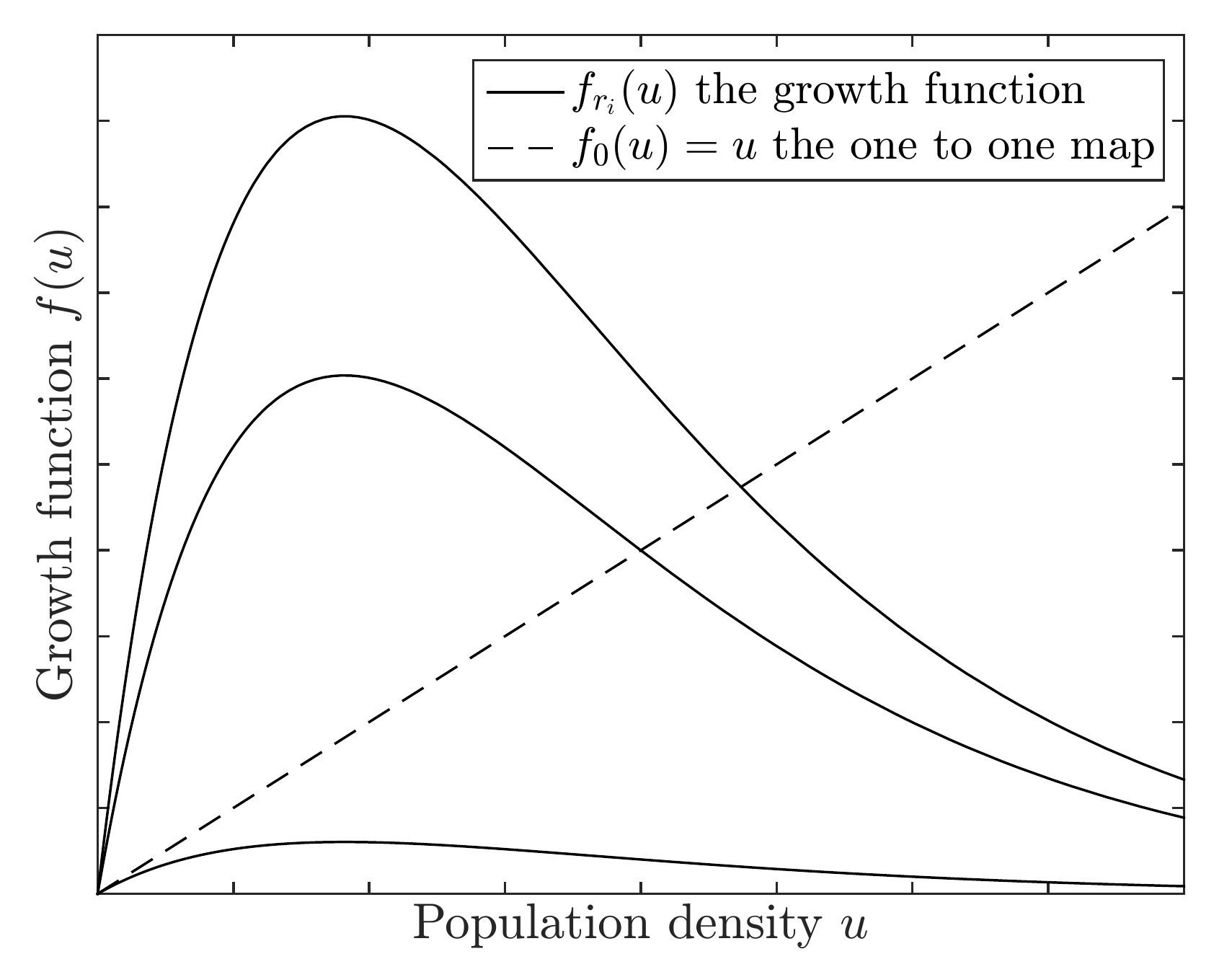}}
\caption{Different realisations for the growth function $f_r$ depending on the random process $r$ for compensatory (a) and over-compensatory (b) growth. One notices that for some realisation $f'_r(0)<1$ and there is no positive fixed point for this type of function $f_r$.}\label{frand}
\end{center}
\end{figure}

We incorporate stochasticity into our modelling framework in two ways, first with respect to the shifting speed for suitable habitat and second with respect to the growth rate within the suitable habitat.   The model for the shifting speed gives the center of the suitable habitat as $s_t=ct+\sigma_t$, where $\sigma_t$ is a random variable. Here we assume that $c$ is unknown but fixed, depending on the scenario considered for the severity of global warming ($c\in\{c_1,\dots,c_n\}$) (Figure \ref{fig:CSigma}).  Variability of the growth dynamics at each generation, for example due to weather conditions or extreme climate events, is included through stochasticity in the growth function $f(u)$. Our approach is to incorporate the randomness into the geometric growth rate $r_t\in\{r_1,\dots,r_n\}$ (Figure \ref{frand}) and so in any given year the growth function is given by $f_{r_t}(u)$ where $r_t=f_{r_t}'(0)$. 

Our goal is to mathematically analyze the effect of climate change and shifting range on population dynamics for a species that grows and then disperses at each generation, taking into account stochasticity induced by environmental variability and climate change. In Section \ref{sec:model} we derive the model and state model assumptions. 
In Section \ref{sec:pers} we first define a persistence condition, derived from the papers of \citet{Hardinetal1} and \citet{JJL}. We then detail computation of the persistence criterion and highlight the link between persistence of the population and principal eigenvalue of the operator linearised about 0 for the case where the shifting speed is not random ($\sigma_t\equiv0$). We also prove that, when the dispersal kernel is Gaussian and the speed is not random, there exists a critical shifting speed characterising persistence of the population. Speeds above this critical value will drive the population to extinction, while speeds below will allow the population to persist. In section \ref{sec:example} we apply the theory to an example in butterfly population subject to changing temperatures in Canada \citep{Lerouxetal13} and use numerical simulations to investigate the dependence of the critical shifting speed on the variance of the dispersal Gaussian kernel. We also compute the persistence criterion as a function of the variance of the dispersal kernel for Gaussian and Laplace kernel. Lastly we numerically investigate the effect of the variability in the yearly shift and in the growth rate on the persistence of the population. In the Appendix we draw on classical theory of spreading speeds for stochastic integrodifference equations to aid development of a heuristic link between the speed of the stochastic wave in the homogeneous framework and the critical domain size and persistence condition.

\section{The model}\label{sec:model}

In this section we derive the mechanistic model used to study population persistence facing global warming and habitat shifts. We explain the different assumptions made for each component throughout the derivation of the model and conclude by explaining how it results in the problem in a moving environment.

As already stated in the introduction, we use the theory of integrodifference equations to model the temporal dynamics of the density of the population $u$, as introduced by \citet{KS86}. In the classical homogeneous case, the model is as given in equation (\ref{eq:IDE}) with $\Omega$ given by $\R$ and with $u_0$ given, nonnegative, compactly supported and bounded.

We make the following assumptions regarding the dispersal kernel $K$
\begin{hyps}\label{hypK}:
\begin{enumerate}[(i)]
\item For all $\xi\in\R$, $\eta\in\R$,
$K(\xi,\eta)=K(\xi-\eta)$.
\item 
$K(x)$ is well defined, continuous, uniformly bounded and positive in $\R$.
\end{enumerate}
\end{hyps}
The first hypothesis means that $K$ takes the form of a difference kernel and depends only on the signed distance between $\xi$ and $\eta$. The second holds for typical kernel, such as the Laplace and Gaussian (Figure \ref{fig:kernel}).

To include the effect of climate change on range distribution in this model, we assume that the suitability of the environment is heterogeneous in the sense that % for all $\xi\in\R$, $t\in\N$,
\be u_{t+1}(\xi)=\int_\R K(\xi-\eta)g_t(\eta)f(u_t(\eta))d\eta.\ee
where the function $(t,\eta)\mapsto g_t(\eta)$ stands for the suitability of the environment at generation $t$ and location $\eta$. 

We make the following assumptions for the suitability function $g_t$:
\begin{hyps}\label{hyp:Q}:
\begin{enumerate}[(i)]
\item Denoting $s_t$ as the reference point on a suitable habitat, we assume that $g_t(\eta)=g_0(\eta-s_t),$ 
\item $g_0(x)$ is compactly supported, nonnegative, bounded by 1 and is non trivial in $\R$. 
\end{enumerate}
\end{hyps}
Biologically these assumptions mean that the suitable environment has a constant profile $g_0$ that is shifted by $s_t$ at generation $t$. 

The model then becomes
\be u_{t+1}(\xi)=\int_\R K(\xi-\eta)g_0(\eta-s_t)f(u_t(\eta))d\eta.\ee
We denote $\Omega_0$ as the support of $g_0$, i.e $\Omega_0:=\{x\in\R,\:g_0(x)>0\}$.
Notice that only the population located in $\Omega_0+s_t:=\{x\in\R\:| \: x=x'+s_t,\:  x'\in\Omega_0\}$ contributes to the growth from generation $t$ to generation $t+1$.  
To introduce environmental stochasticity in our model, we assume that $(s_t)_{t\in\N}$ is a random process and $f(u)=f_{r_t}(u),$ with $(f_{r_t})_{t\in\N}$ a sequence of random functions. We interpret $r_t$ as the geometric growth rate of the population at low density. These two forms of environmental stochasticity emphasise the dependence of the range shift and the growth rate on the strength of the climate change. 
Moreover we can be more precise about the form of the shift variable $s_t$. Indeed, we assume that for all $t\in\N$, $s_t=ct+\sigma_t,$ where $c>0$ is a constant representing the asymptotic shifting speed and $\sigma_t$ is a random variable representing the environmental stochasticity in the shift from one year to the next. In the introduction we stated that the asymptotic shifting speed itself may be uncertain. However from now on we consider it to be a fixed constant $c$ and study the problem of persistence of the population for different possible values of the asymptotic shifting speed $c$.

Denoting by $(\alpha_t)_t=(\sigma_t,r_t)_t$ and $\mathcal{S}$ the set of possible outcomes for $\alpha$ at each generation, we assume that the elements $(\alpha_t)_t$ are independent, identically distributed and bounded by appropriate values, namely:
\begin{hyps}\label{alpha}:
\begin{enumerate}[(i)]
\item $(\alpha_t)_t=(\sigma_t,r_t)_t$ is a sequence of independent, identically distributed random variables, with distribution $\mathcal{P}_\alpha$ %\label{alpha1}\\
\item There exists $\underline{\sigma}<0<\overline{\sigma}$ such that for all $t\in\N,\: \underline{\sigma}\leq \sigma_t<\overline{\sigma}$ with probability 1%}.\label{alpha2}\\
\item There exists $\overline{r}>\underline{r}>0$ such that for all $t\in\N,\: \underline{r}\leq r_t<\overline{r}$ with probability 1%}.\label{alpha3}
\end{enumerate}\end{hyps}
We consider a self-regulating population, with negative density-dependence and make the following assumptions on the growth function $f$:
\begin{hyps}\label{hypf}
For any $r$ such that $\alpha\in\mathcal{S}$ 
\begin{enumerate}[(i)]%\label{hyp:f1}
\item $f_r :\R\to[0,+\infty)$ is continuous, with $f_r(u)=0$ for all $u\leq0$,
\item There exists a constant $m>0$ such that, for all $r$,
\begin{enumerate}[a.]
%\begin{gather}\label{hyp:f2}
%\be\label{f2-1}
\item $u\in\R^+\mapsto f_r(u)$ is nondecreasing,%\ee
%\be\label{f2-2}
\item $0< f_r(u)\leq m$ for all positive continuous function $u$,%\ee
%\be\label{f2-3}\text
\item If $u,\:v$ constants such that $0<v<u$ then $f_r(u)v<f_r(v)u$,%},\ee
%\be\label{f2-4}
\item $u\in\R\mapsto f_r(u)$ is right differentiable at  0, uniformly with respect to $\alpha\in\mathcal{S}$.%}.\ee
\end{enumerate}
%\end{gather}
\item We denote $r=f'_r(0)$ as the right derivative of $f_r$ at 0 and assume for now that
\begin{enumerate}[a.]
%\begin{gather}\label{hyp:f3}
%\be\label{f3-1}
\item $\underline{r}\leq f'_r(0)\leq \overline{r}$, from Hypothesis \ref{alpha}(iii)%\ee
%\be\label{f3-2}\text{ 
\item $\inf\{f_r(b),\: \alpha\in\mathcal{S}\}>0$, with $b:=m \sup_{x\in\R}{K(x)}\int_{\Omega_0} g_0(y)dy$.%}.\ee
\end{enumerate}
\end{enumerate}
\end{hyps}
Hypothesis \ref{hypf}(i) means that the population does not grow when no individuals are present in the environment. We also assume that the growth of the population is bounded (Hypothesis \ref{hypf}(ii)b.) and  consider a population not subject to an Allee effect, that is, the geometric growth rate is decreasing with the density of the population (Hypothesis \ref{hypf}(ii)c.). Hypothesis \ref{hypf}(ii)d. holds for typical growth functions as the ones illustrated by Figure \ref{fgrowth}.  Notice that because of Hypothesis \ref{hypf}(ii)a., we do not yet consider models with over-compensatory competition here. Indeed we will need this monotonicity assumption to prove large time convergence of the solution of our problem. Nevertheless, we can extend our result to the case with over-compensation where $f$ is not assumed to be nondecreasing and so Hypothesis \ref{hypf}(ii)a. no longer holds. More details are given in Corollary \ref{cor:overcomp}, stated in the next section. Hypothesis \ref{hypf}(iii)a. assumes that the geometric growth rate at zero density is bounded from above and below by some positive constants, while Hypothesis \ref{hypf}(iii)b. implies that the growth term at the maximal density is positive, for all the possible environments. 

Finally, denoting by $F_{\alpha_t}(u):=\int_\Omega K(x-y+c)g_0(y-\sigma_t)f_{r_t}(u(y))dy$, we make a last assumption: 
\begin{hyp}\label{alphastar}
There exists $\alpha^*\in\mathcal{S}$ such that 
%\be\label{alphastar}
\be F_\alpha(u)\leq F_{\alpha^*}(u)\ee
%\ee
for all $\alpha\in\mathcal{S}$, $u$ nonnegative continuous function.
\end{hyp}
This assumption means that there exists an environment $\alpha^*$ that has the better outcome than all other environments in terms of population growth.

The general problem becomes
 \be\label{IDEnonmoving}
u_{t+1}(\xi)=\int_\R K(\xi-\eta)g_0(\eta-s_t)f_{r_t}(u_t(\eta))d\eta.
\ee

We are interested in the large time behaviour of the density $u$. Using a similar approach to the one of \citet{ZK1}, we would like to study the problem in the shifted environment, so as to track the population. Considering the moving variables $\xi_{t+1}:=\xi-s_{t+1}$, $\eta_t:=\eta-s_t$ and letting $\overline{u}_t(\eta_t):=u_t(\eta+s_t)$ be the associated density in the moving frame (where the reference point is the centre of the suitable environment that is shifted by $s_t$ at time $t$), we obtain a problem where the space variables $\xi_{t+1}$ and $\eta_t$ are  also random variables. To simplify the analysis, we consider the problem in the ``asymptotic" moving frame. That is we consider $x=\xi-c(t+1)$ and $y=\eta-ct$ and using Hypothesis \ref{hypK}(i), equation \eqref{IDEnonmoving} can be written as
\be u_{t+1}(x+c(t+1))=\int_\R K(x-y+c)g_0(y-\sigma_t)f_{r_t}(u_t(y+ct))dy.\ee
Letting $\overline{u}_t(y):=u_t(y+ct)$ be the associated density in the moving frame, we have
\be\overline{u}_{t+1}(x)=\int_\R K(x-y+c)g_0(y-\sigma_t)f_{r_t}(\overline{u}_t(y))dy.\ee
Dropping the bar, we obtain the following problem, in the moving frame,
\begin{equation}\label{IDEmoving1}
u_{t+1}(x)=\int_\R K(x-y+c)g_0(y-\sigma_t)f_{r_t}(u_t(y))dy,\\
\end{equation}
where $u_0$ is given as a nonnegative, non trivial and bounded function. 
Moreover, as stated in Hypothesis \ref{alpha}(ii), for all $t\in\N$, $\sigma_t\in(\underline{\sigma},\overline{\sigma})$ and thus defining 
\be\label{Omega}\Omega:=(\inf \Omega_0+\underline{\sigma},\sup \Omega_0+\overline{\sigma}),\ee
we only have to study equation \eqref{IDEmoving1} for all $x\in\Omega$, for all $t\in\N$. 
The problem in the moving frame becomes
 \begin{equation}\label{IDEmoving}
u_{t+1}(x)=\int_\Omega K(x-y+c)g_0(y-\sigma_t)f_{r_t}(u_t(y))dy,\\
\end{equation}
where $u_0$ is given as a nonnegative, non trivial and bounded function.  Notice that this problem is now defined on a compact set $\Omega\subset\R$.

Many results already exist for the deterministic version of \eqref{IDEmoving}. When $\sigma_t\equiv 0$ and $f$ is deterministic, depending only on $u$, it has been shown that, if $g_0(y)=\mathds{1}_{[-\frac{L}{2}, \frac{L}{2}]}$, the magnitude of the largest eigenvalue of the linearised operator around 0 determines the stability of the trivial solution and thus the persistence of the population (see \citet{ZK1, ZK2} for the one dimensional problem, \citet{PK15} for the two dimensional problem). \citet{HZHK} study a similar deterministic integrodifference equation modified to include age and stage structured population and show that if the magnitude of the principal eigenvalue of the linearised problem around 0 exceeds one then 0 becomes stable and the population goes extinct. Zhou also studied the associated deterministic problem \eqref{IDEmoving1}, for more general deterministic functions $f$ and $g_0$ \citep[][chapter 4]{ZhouThesis}. 

In this paper we are interested in similar questions about persistence of the population, but now for a shifting environment that moves at a random speed and a population that reproduces at a growth rate, chosen randomly at each generation.

Mathematically the assumptions on $K$, $f$ and $F$ that we make are similar than those made by \citet{JJL} in their study of the question of persistence of a population in temporally varying river environment. 
They assumed that the dispersal kernel and growth terms are randomly distributed at each time step but also assumed asymmetric dispersal kernels to take the effects of water flow into account. They used the theory of \citet{Hardinetal1} to develop a persistence criterion for the general model
\be u_{t+1}(x)=F_{\alpha_t}(u_t(x)),\ee
where $t\in\N$, $x\in\Omega\subset\R^n$ and $(\alpha_t)_t$ are independently identically distributed random variables and related this persistence criterion to the long-term growth rate. In what follows, we use similar theory to study our integrodifference equation \eqref{IDEmoving} where $\alpha_t=(\sigma_t,r_t)_t$. 

\section{Persistence condition for random environments with shifted kernel}\label{sec:pers}
In this section we derive a criterion that separates persistence from extinction for a population facing random environments and climate change. We highlight the dependence of this criterion on the asymptotic geometric growth rate at low density on one hand and on the shifted dispersal success function on the other hand. We then develop a connection between persistence and the magnitude of the principal eigenvalue of a linearised operator for the case when only the growth rate is stochastic. This allows us to conclude on the existence of a critical shifting speed for the case of Gaussian dispersal kernel.

\subsection{Conditions for persistence}
We first describe the derivation of the deterministic criterion characterising persistence from extinction of the population.

Define $(u_t)_t$ to be the solution of problem \eqref{IDEmoving}. Then $(u_t)_t$, the sequence of population density at each generation $t\in\N$, is a random process and for each $t$ positive, $u_t$ takes values in $C_+(\Omega)$, the set of continuous, nonnegative function defined on $\Omega$. Using the same notation as \citet{Hardinetal1}, we can rewrite equation \eqref{IDEmoving} as follows
\be u_{t+1}(x)=F_{\alpha_t}[u_t](x),\ee
for all $x\in\Omega$.

We have the following theorem about the large time behaviour of the solution
\begin{theorem}\label{th:convergence}
Assume that $K$, $g_0$, $f_r$ and $(\alpha_t)_t$ satisfy all the assumptions stated in Section \ref{sec:model} (Hypotheses \ref{hypK} - \ref{alphastar}), let $(u_t)_t$ be the solution of problem \eqref{IDEmoving}, with a bounded initial condition $u_0$. Then $u_t$ converges in distribution to a random variable $u^*$ as time goes to infinity, independently of the initial condition $u_0$, and $u^*$ is a stationary solution of \eqref{IDEmoving}, in the sense that
\be u^*(x)=\int_\Omega K(x-y+c)g_0(y-\sigma^*)f_{r^*}(u^*(y))dy,\ee
with $(\sigma^*,r^*)$ a random variable taking its values in $\mathcal{S}$ with distribution $\mathcal{P}_\alpha$ (Hypothesis \ref{alpha}(i)).
Denoting by $\mu^*$ the stationary distribution associated with $u^*$ and $\mu^*(\{0\})$ the probability that $u^*\equiv0$, we have that $\mu^*(\{0\})=0$ or $\mu^*(\{0\})=1$
\end{theorem}
The proof of this theorem follows from \citet[Theorem 4.2] {Hardinetal1} and is detailed in the Appendix \ref{sec:proofhardinetal}. One can then be more precise about the distribution $\mu^*$. Define 
\be\label{Lambda}\Lambda_t:=\left(\int_\Omega \tilde{u}_t(x)dx\right)^{1/t},\ee
where $(\tilde{u}_t)_t$ is the solution of the linearised problem around 0, i.e for all $x\in\Omega$
\be\label{eq:utilde}\tilde{u}_{t+1}(x)=\mathcal{L}_{\alpha_t} \tilde{u}_t(x):=\int_\Omega K(x-y+c)g_0(y-\sigma_t)r_t\tilde{u}_t(y)dy.\ee
The metric $\Lambda_t$ can be interpreted as the growth rate of the linearised operator up to time $t$ and its limit, representing the asymptotic growth rate of the linearised operator.
We have the following theorem
\begin{theorem}\label{th:persistence}
Let $\Lambda_t$ be defined as in \eqref{Lambda}, then 
$$\underset{t\to+\infty}{\lim}\Lambda_t=\Lambda\in[0,+\infty) \text{, with probability 1.}$$
And,
\begin{itemize}
\item If $\Lambda<1$, the population will go extinct, in the sense that $\mu^*(\{0\})=1$,
\item If $\Lambda>1$, the population will persist, in the sense that $\mu^*(\{0\})=0$.
\end{itemize}
\end{theorem}
The proof of this theorem follows from \citet[Theorem 5.3]{Hardinetal1} and \citet[Theorem 2]{JJL}, and it is explained in the Appendix \ref{sec:proofhardinetal}.
From these two theorems one can then deduce a corollary, for the case where the growth function $f$ is not assumed to be monotonic anymore.
\begin{corollary}\label{cor:overcomp}
Let $K$, $f$, $g_0$ and $(\alpha_t)_t$ satisfy the assumptions from Section \ref{sec:model} except assumption \ref{hypf}(ii)a, that is, $f_r$ is not assumed to be nondecreasing. 
Let $(u_t)_t$ be the solution of \eqref{IDEmoving}, with a bounded initial condition $u_0$. Let $\Lambda_t$ be defined as in \eqref{Lambda}, then 
$$\underset{t\to+\infty}{\lim}\Lambda_t=\Lambda\in[0,+\infty) \text{, with probability 1.}$$
Also,
\begin{itemize}
\item If $\Lambda<1$, the population will go extinct, in the sense that $\underset{t\to+\infty}{\lim}\: u_t(x)=0$ for all $x\in\Omega$ with probability 1,
\item If $\Lambda>1$, the population will persist, in the sense that $\underset{t\to+\infty}{\liminf}\: \underset{x\in\Omega}{u_t(x)}>0$ with probability 1.
\end{itemize}
\end{corollary}
\begin{proof} %Indeed considering $f_r$ that are not necessarily increasing. 
First note that
\begin{align}
||F_\alpha(u)||_\infty&=\underset{x\in\Omega}{\max}\int_\Omega K(x-y+c)g_0(y-\sigma)f_r(u(y))dy\\
&< m \cdot \sup_{x\in\R}\:K(x)\cdot\int_{\Omega_0} g_0(y)dy=b,
\end{align}
where $m$ is defined in Hypothesis \ref{hypf}(ii)b.
This implies that 
\be u_t(x)\leq b,\ee
for all $t\in\N^*$, for all $x\in\Omega$. Then define lower and upper nondecreasing functions that satisfies Hypotheses \ref{hypf}
$$u\in C_+(\R)\mapsto \underline{f}_r(u),\: u\in C_+(\R)\mapsto \overline{f}_r(u),$$
so their slopes match that of $f_r$ at zero
\be\label{fprimesame}\underline{f}'_r(0)=f'_r(0)=\overline{f}'_r(0)=r,\ee
and they satisfy the inequalities
\be\label{fordered}0<\underline{f}_r(u)\leq f_r(u)\leq \overline{f}_r(u)<m\ee
for all $u\in(0,b)$. For instance, we can choose the nondecreasing function %such that $\forall u\geq0$,
\begin{equation}\label{fminmax}
\underline{f}_r(u)=\min_{u\leq v\leq b}f_r(v) \text{ and } \overline{f}_r(u)=\max_{0\leq v\leq u}f_r(v),\end{equation} 
%and $\underline{f}_r(u)=\overline{f}_r(u)=0$ for all $u<0$. 
One can easily prove, using the definition of $\underline{f}_r$, $\overline{f}_r$ and the properties of $f_r$, that $\underline{f}_r$ and $\overline{f}_r$ as defined in \eqref{fminmax} satisfy Hypotheses \ref{hypf}(i),(ii) and (iii) and are such that \eqref{fprimesame}-\eqref{fordered} are satisfied. 

Then denoting $(v_t)_t$, respectively $(w_t)_t$, as the solution of problem \eqref{IDEmoving}, with $\underline{f}_r$ instead of $f_r$, respectively $\overline{f}_r$ instead of $f_r$, such that $v_0\equiv u_0\equiv w_0$, we have
\be\label{eq:vuw}v_t(x)\leq u_t(x)\leq w_t(x)\ee
for all $t\in\N$, for all $x\in\Omega$. 
Applying the previous theorems (Theorem \ref{th:convergence} and Theorem \ref{th:persistence}) we know that 
\begin{itemize}
\item If $\Lambda>1$ (where $\Lambda$ is defined in Theorem \ref{th:persistence}), then the population associated with the process $(v_t)_t$ persists in the sense that $\mu^*_v(\{0\})=P(\underset{t\to+\infty}{\lim}v_t\equiv0)=0$, which implies, using \eqref{eq:vuw} that $\underset{t\to+\infty}{\lim}\: \underset{x\in\Omega}{\inf}u_t(x)>0$ with probability 1,
\item If $\Lambda<1$, then the population associated with the process $(w_t)_t$ goes extinct in the sense that $\mu^*_w(\{0\})=P(\underset{t\to+\infty}{\lim}w_t\equiv0)=1$, which implies, using again \eqref{eq:vuw}, that $\underset{t\to+\infty}{\lim}\: u_t(x)=0$, for all $x\in\Omega$, with probability 1.
\end{itemize}
This completes the proof of Corollary \ref{cor:overcomp}.
\end{proof}
This means that the results about persistence and extinction extend to the case of over-compensatory growth as illustrated in Figure \ref{fgrowth}, even though we cannot draw conclusions about the large time convergence of the solution. 

Note that the persistence criterion as given by Theorem \ref{th:persistence} or Corollary \ref{cor:overcomp} is difficult to analyse further. It requires the calculation of the asymptotic growth rate of the linearised operator under random environmental conditions (equation \eqref{Lambda}). However we can proceed heuristically regarding a necessary condition for persistence of the population by comparing the asymptotic shifting speed $c$ with the asymptotic spreading speed of the population in a homogeneous environment. Consider the population dynamics in a homogeneous environment ($g_0\equiv1$). Equation \eqref{eq:utilde} becomes
\be\label{eq:utildehomo}
\tilde{u}_{t+1}(x)=\int_\R K(x-y+c)r_t\tilde{u}_t(y)dy.
\ee
If we consider exponentially decaying solution of the form $\tilde{u}_t\propto h_te^{-sx}$, $s>0$, substitution into \eqref{eq:utildehomo} yields
\be h_{t+1}=M(s)e^{-sc}h_tr_t\ee
where $M$ is the moment generating function of $K$
\be\label{MomGenFunc} M(s)=\int_\R e^{sz}K(z)dz.\ee
The pointwise growth rate of the solution is given asymptotically by 
\begin{align} \underset{t\to+\infty}{\lim}\frac{1}{t} \ln\left(\frac{h_t}{h_0}\right)&=\ln(M(s))+\underset{t\to+\infty}{\lim}\:\frac{1}{t}\sum_{i=0}^{t-1} \ln(r_i)-cs\\
&=E[\ln(r_0M(s))]-cs\end{align}
and is positive if $c<\overline{c}^*$ where 
\be \overline{c}^*=\inf_{s>0}\frac{1}{s} E[\ln(r_0 M(s))].\ee
We interpret $\overline{c}^*$ as the spreading speed of a population undergoing dispersal and stochastic growth in a spatially homogeneous environment (Appendix \ref{a:stochspeed}). In the case where $g_0$ is diminished, so it is not equal to one at all point in space, $c<\overline{c}^*$ may no longer suffice to give a positive pointwise growth rate.
On the other hand, if we consider a bounded, compactly supported initial condition, $\tilde{u}_0$, and thus for all $s>0$ there exists $h_0$ such that $\tilde{u}_0(x)\leq h_0 e^{-sx}$, then the population $\tilde{u}_t$, solution of \eqref{eq:utildehomo} starting with the initial condition $\tilde{u}_0$, 
\be \tilde{u}_{t+1}(x)\leq h_{t+1}e^{-sx}.\ee
The pointwise growth rate of the solution $\tilde{u}_t$ will thus be smaller that the pointwise growth rate of the solution $h_te^{-sx}$ and $c<\overline{c}^*$ is  necessary to give a positive pointwise growth rate of $\tilde{u}_t$.

\subsection{Computation of the persistence criterion $\Lambda$}\label{sec:approxlambda}
In this section, we detail possible computation of the persistence criterion $\Lambda$ and highlight the link between $\Lambda$ and the principal eigenvalue of the linearised operator around zero. 

Using the analysis of the previous section one can determine whether a population can track its favorable environment and thus persist ($\Lambda>1$ in Theorem \ref{th:persistence} and Corollary \ref{cor:overcomp}) or cannot keep pace with the shifting environment and goes extinct ($\Lambda<1$ in Theorem \ref{th:persistence} and Corollary \ref{cor:overcomp}). It would be interesting to understand the dynamics of $\Lambda:=\underset{t\to+\infty}{\lim} \Lambda_t$ in terms of the different parameters of the problem. Using the definition of $\Lambda_t$ \eqref{Lambda}, we have
\begin{align}
\Lambda_1&=\int_\Omega \tilde{u}_1(x)dx\\
&=\int_\Omega\int_{\Omega}K(x-y+c)g_0(y-\sigma_0)r_0u_0(y)dydx
\end{align}
where $r_0$, $\sigma_0$ are the realisations for the geometric growth rate and the yearly shift at generation 0.  Continuing,
\begin{align}
\Lambda_2&=(\int_\Omega \tilde{u}_2(x)dx)^{1/2}\\
&=(\int_\Omega\int_{\Omega}K(x-y_1+c)g_0(y_1-\sigma_1)r_1\tilde{u}_1(y_1)dy_1dx)^{1/2}\\
&=(r_1r_0)^{1/2}(\int_\Omega\int_{\Omega}\int_\Omega K(x-y_1+c)g_0(y_1-\sigma_1)K(y_1-y_2+c)g_0(y_2-\sigma_0)u_0(y_2)dy_2dy_1dx)^{1/2}
\end{align}
where $r_1$, $\sigma_1$ are the realisations for the geometric growth rate and the yearly shift at generation 1. 
In a similar manner, one can then derive the more general formula
\be\label{Lambdat}\Lambda_t=R_t^{1/t}\cdot \kappa_t^{1/t}\ee
where 
\be R_t=\prod_{i=0}^{t-1} r_i\in\R^+,\ee
the geometric mean of the geometric growth rate at zero and 
\begin{align}
\kappa_t&=\underbrace{\int_\Omega\dots\int_\Omega}_{t+1 \text{ terms}}K(x-y_1+c)g_0(y_1-\sigma_{t-1})\cdots K(y_{t-1}-y_{t}+c)g_0(y_{t}-\sigma_0)u_0(y_t)dy_t\dots dy_{1}dx
\end{align}
Considering 
\be \ln(R_t^{1/t})=\frac{1}{t}\sum_{i=0}^{t-1}\ln(r_i),\ee
and using the strong law of large number we have 
$$\ln(R_t^{1/t})\to E[\ln(r_0)]\in\R \text{ as } t\to+\infty$$
with probability 1 and thus obtain
\be\label{Rbar}R_t^{1/t}\to e^{E[\ln(r_0)]}=\overline{R}\in(0,+\infty) \text{ as } t\to+\infty,\ee
with probability 1. Combining equations \eqref{Lambdat} and \eqref{Rbar}, the persistence criterion $\Lambda$ becomes
\be\label{Lambdar0K}\Lambda=e^{E[\ln(r_0)]}\cdot \lim_{t\to+\infty}\kappa_t^{1/t}.\ee
This formula for $\Lambda$ highlights the dependence of the persistence criterion on the different parameters of the problem. Notice that the distribution of the growth rate at zero affects the first part whereas the distribution of the yearly shifts affects the second part of the formula. We further analyse the effect of the variation in $r$ or the variation in $\sigma$ using numerical simulation in Section \ref{sec:example}.

\subsubsection*{Persistence criterion $\Lambda$ for deterministic shifting speed}
We turn our attention to further analysis of the persistence criterion $\Lambda$, as given in \eqref{Lambdar0K}, when we assume that the shifting speed is not random. That is, we assume that $\sigma\equiv0$ and thus randomness comes only from the growth term. We still consider a population that sees its favorable environment shifted at a speed $c$, but this speed is not assumed to have random variation from one generation to the next. In this case $\Omega=\Omega_0$ and 
\be\Lambda_t=R_t^{1/t}\cdot \left(\int_{\Omega_0} \cdots\int_{\Omega_0} K(x-y_1+c)g_0(y_1) \dots K(y_{t-1}-y_{t}+c)g_0(y_{t})u_0(y_t)dy_t\dots dy_{1}dx\right)^{1/t}.\ee
Define $\mathcal{K}_c$ the linear operator such that 
\be\label{Kc}\mathcal{K}_c[u](x)=\int_{\Omega_0}K(x-y+c)g_0(y)u(y)dy.\ee
As $K$ is positive and $\Omega_0$ is compact, the operator $\mathcal{K}_c$ is compact \citep{K64} and strongly positive, that is for any function $u\geq0$, there exists $t\in\N$ such that $\mathcal{K}_c^t[u](x):=\mathcal{K}_c\bigg[\mathcal{K}_c\Big[\dots[\mathcal{K}_c[u]\dots\Big]\bigg](x)>0$ for all $x\in\Omega_0$. Then applying the Krein-Rutman theorem, it follows that this operator possesses a principal eigenvalue $\lambda_c>0$ such that $|\lambda|<\lambda_c$ for all other eigenvalues and $\lambda_c$ is the only eigenvalue associated with a positive principal eigenfunction $\phi_c$. That is $\lambda_c>0$ and $\phi_c>0$ satisfy
\be\lambda_c \phi_c(x)=\int_{\Omega_0}K(x-y+c)g_0(y)\phi_c(y)dy\ee
for all $x\in\Omega_0$. One can always normalise $\phi_c$ such that $\int_{\Omega_0}\phi_c(y)dy=1$.
Now choosing the initial condition for equation \eqref{IDEmoving} so that $u_0\equiv\phi_c$, we have
\be\Lambda_t=R_t^{1/t}\cdot \lambda_c\ee
so that 
\be\Lambda=\overline{R}\cdot\lambda_c,\ee
with $\overline{R}= e^{E[\ln(r_0)]}$ defined in \eqref{Rbar}. 

One has to approximate the principal eigenvalue $\lambda_c$ to be able to conclude about the persistence of the population. One can refer to the paper by \citet{KP15} for the description and implementation of different numerical or analytical methods to compute the principal eigenvalue of a linear operator.

In the specific case of Gaussian kernel with variance $(\sigma^K)^2$, that is
\be K(x)=\frac{1}{\sqrt{2\pi (\sigma^K)^2}}e^{-\frac{x^2}{2(\sigma^K)^2}},\ee
the principal eigenvalue of the shifted linear operator $\mathcal{K}_c$ can be expressed as a decreasing function of $c$. Indeed, we know that $\lambda_c$ and $\phi_c$ satisfy for all $x\in\Omega_0$
\begin{align}
\lambda_c \phi_c(x)&=\int_{\Omega_0}\frac{1}{\sqrt{2\pi (\sigma^K)^2}}e^{-\frac{(x-y+c)^2}{2(\sigma^K)^2}}g_0(y)\phi_c(y)dy\\
&=\int_{\Omega_0}\frac{1}{\sqrt{2\pi (\sigma^K)^2}}e^{-\frac{(x-y)^2}{2(\sigma^K)^2}}e^{-\frac{c^2}{2(\sigma^K)^2}}e^{-\frac{2c(x-y)}{2(\sigma^K)^2}}g_0(y)\phi_c(y)dy\\
\end{align}
The last equality is equivalent to 
\be (e^\frac{c^2}{2(\sigma^K)^2}\lambda_c)\cdot(e^{\frac{cx}{(\sigma^K)^2}}\phi_c(x))=\int_{\Omega_0}\frac{1}{\sqrt{2\pi (\sigma^K)^2}}e^{-\frac{(x-y)^2}{2(\sigma^K)^2}}g_0(y)(e^{\frac{cy}{(\sigma^K)^2}}\phi_c(y))dy.\ee
Thus one has that 
\be\lambda_c=e^{-\frac{c^2}{2(\sigma^K)^2}}\lambda_0,\ee
where $\lambda_0$ is the principal eigenvalue of the linear operator $\mathcal{K}_0$ (defined by \eqref{Kc} when $c=0$) and %. As 
%$$E[\Lambda_t]=E\left[r_0^{1/t}\right]^t\cdot \lambda_c,$$
%we have that 
\be\label{LambdaGauss}\Lambda=e^{-\frac{c^2}{2(\sigma^K)^2}}\lambda_0\cdot\overline{R}\ee
is decreasing with $c>0$.
Assuming that the parameters of the problem are chosen so that when $c=0$ the population persists, that is 
\be \overline{R}\cdot\lambda_0>1,\ee
there exists a critical value $c^*>0$ such that for all $c<c^*$ the population persists and when $c>c^*$ the population goes extinct (in the sense of Theorem \ref{th:persistence} or Corollary \ref{cor:overcomp}). The critical value $c^*$ is such that 
\be e^{-\frac{(c^*)^2}{2(\sigma^K)^2}}\lambda_0\cdot\overline{R}=1\ee
that is
\be c^*=\sqrt{2(\sigma^K)^2\left(\ln(\lambda_0)+E[\ln(r_0)]\right)}.\ee
We can thus state the following proposition
\begin{proposition}\label{CritSpeedGauss}
Let $K$, $f$ and $F$ satisfy the assumption from Section \ref{sec:model}. Assume also that $\sigma_t\equiv 0$ and $K$ is a Gaussian kernel with variance $(\sigma^K)^2$, that is
\be K(x)=\frac{1}{\sqrt{2\pi (\sigma^K)^2}}e^{-\frac{x^2}{2(\sigma^K)^2}}.\ee
There exists $c^*\geq0$ such that
\begin{itemize}
\item for all $c<c^*$, the population persists in the sense of Theorem \ref{th:persistence} (or Corollary \ref{cor:overcomp} when the growth map $f$ is not necessarily monotonic),
\item for all $c>c^*$, the population goes extinct in the sense of Theorem \ref{th:persistence} (or Corollary \ref{cor:overcomp} when the growth map $f$ is not necessarily monotonic).
\end{itemize}
Moreover, if $\overline{R}$, defined in \eqref{Rbar}, and $\lambda_0$ the principal eigenvalue of $\mathcal{K}_0$ defined in \eqref{Kc} are such that $\overline{R}\cdot\lambda_0>1$, the critical shifting speed is given by
\be\label{eq:cstar}c^*=\sqrt{2(\sigma^K)^2\left(\ln(\lambda_0)+E[\ln(r_0)]\right)}>0.\ee
\end{proposition}

The critical shifting speed for persistence thus increases with the expected geometric growth rate. Heuristically one would imagine that when $\sigma^K$ is small enough relative to $\Omega_0$ then $c^*$ should increase with the variance of the dispersal kernel, as the population becomes more mobile. On the other hand when $\sigma^K$ becomes too large relative to $\Omega_0$ then the population would disperse outside its favorable environment and then $c^*$ should decrease with the variance of the dispersal kernel.

\begin{figure}\begin{center}
\subfloat[]{\includegraphics[scale=0.35]{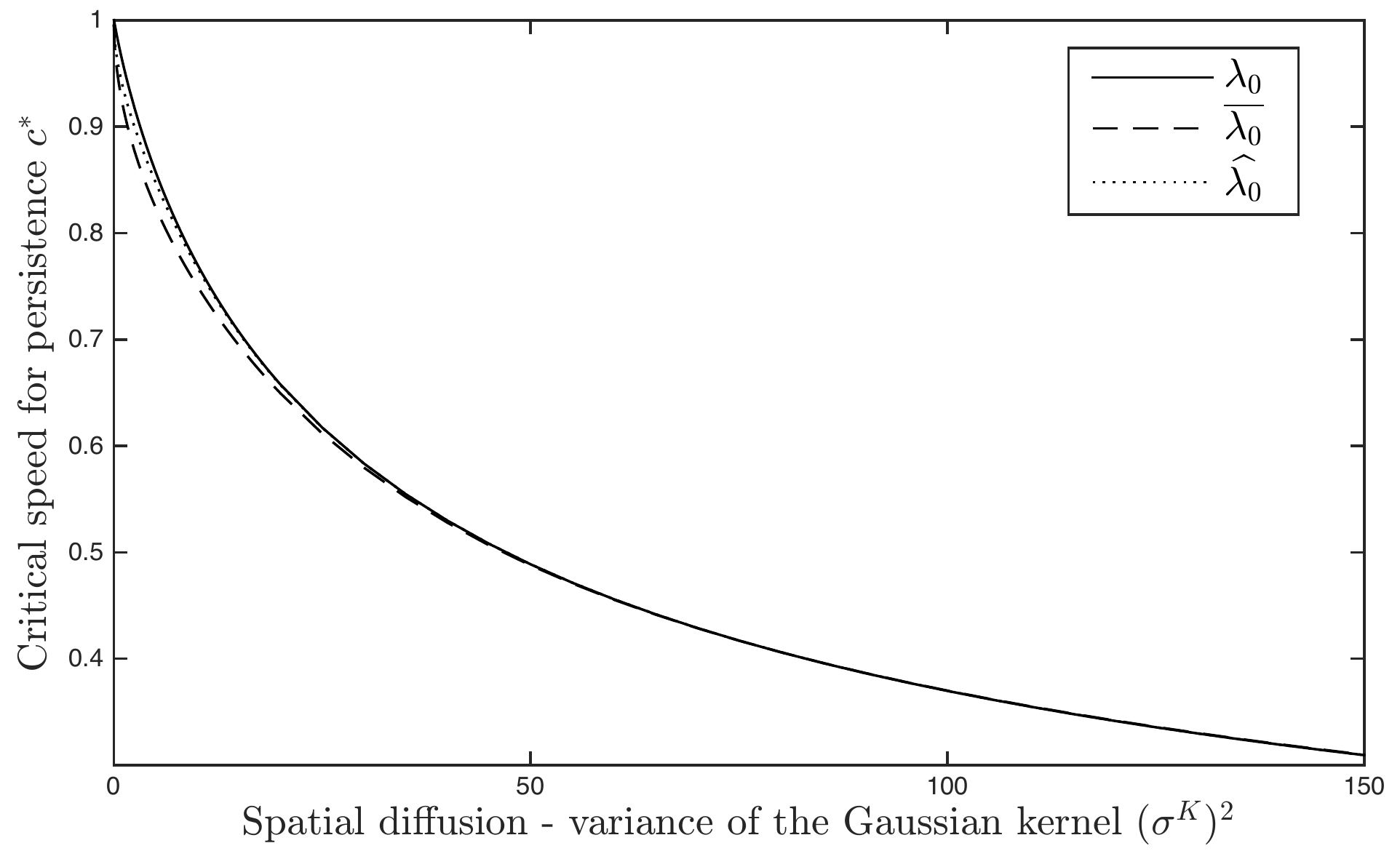}}\quad
\subfloat[]{\includegraphics[scale=0.35]{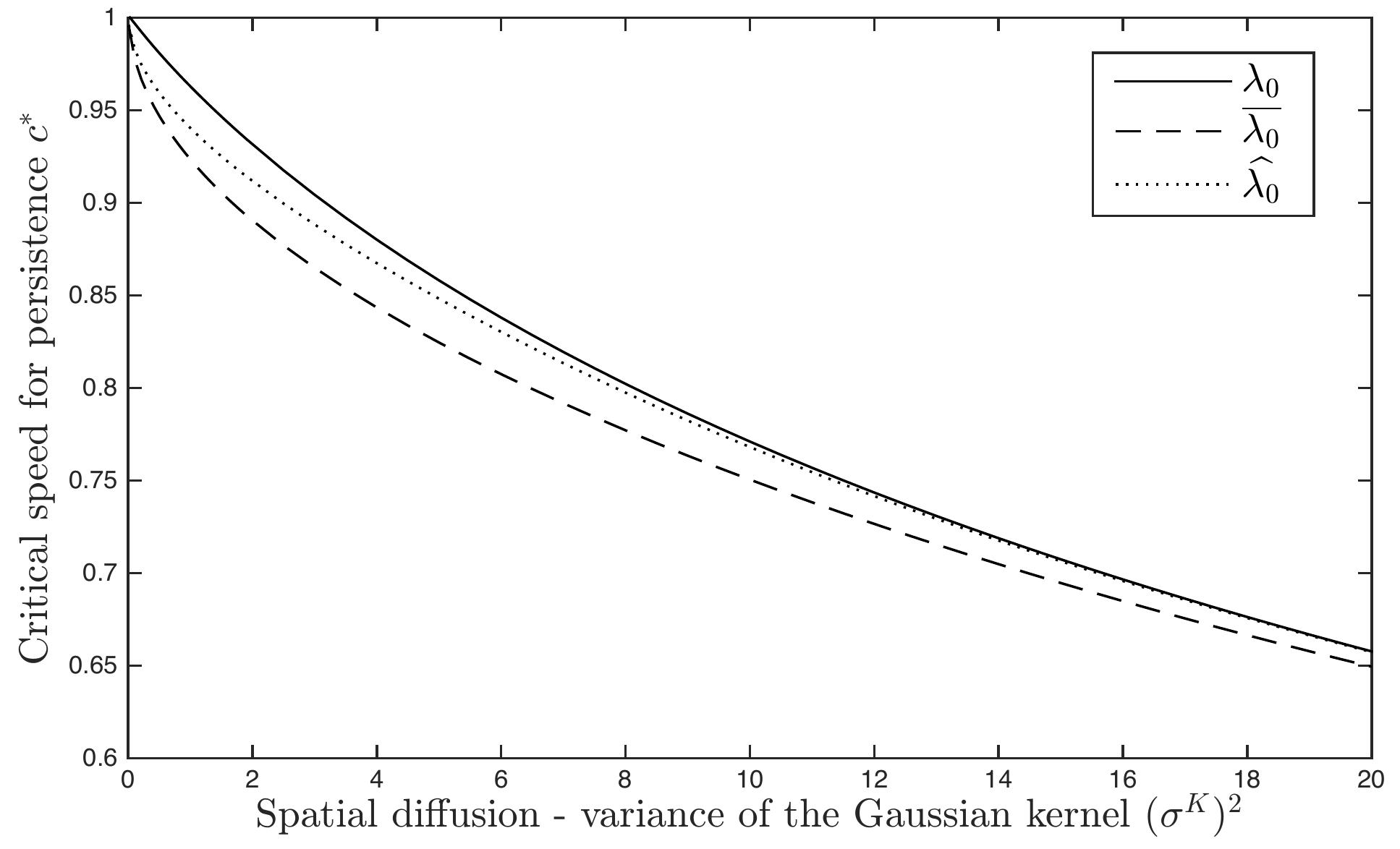}}
\caption{The principal eigenvalue $\lambda_0$, its dispersal success approximation $\overline{\lambda_0}$ and its modified dispersal success approximation $\widehat{\lambda_0}$ as a function of $(\sigma^K)^2$, the variance of the dispersal kernel, when $K$ is a gaussian kernel and $|\Omega_0|=10$. On the right panel (a), $(\sigma^K)^2$ spans from 0.01 to 150. The left panel (b) shows a zoom of the right panel for small values of $(\sigma^K)^2$, where the error between the different approximation is the largest.}\label{ppleeigenDispSuccAppr}
\end{center}
\end{figure}

Assuming that
$$g_0(y)=\mathds{1}_{\Omega_0}(y)=\begin{cases} 1 &\text{if }y\in\Omega_0,\\
									0 &\text{otherwise.}\end{cases}$$
one can also approximate the principal eigenvalue $\lambda_0$ using the dispersal success approximation 
\be\overline{\lambda_0}=\frac{1}{|\Omega_0|}\int_{\Omega_0}s^0(y)dy\ee
as defined in \eqref{eq:DispSuccAppr} or the modified dispersal success approximation 
\be\widehat{\lambda_0}=\frac{1}{\Omega_0}\int_{\Omega_0}\left(\frac{s^0(y)}{\overline{\lambda_0}}\right)s^0(y)dy\ee
as defined in \eqref{eq:ModDispSuccAppr}, with $s^0(y)=\int_{\Omega_0}K(x-y)dx$ the dispersal success function when $c=0$. 
We obtained similar conclusion than \citet{RBM15} when comparing the principal eigenvalue $\lambda_0$, its dispersal success approximation $\overline{\lambda_0}$ and its modified dispersal success approximation $\widehat{\lambda_0}$ (Figure \ref{ppleeigenDispSuccAppr}), observing that the modified dispersal approximation gives a more accurate estimate of the principal eigenvalue $\lambda_0$ than the regular dispersal success approximation.

In this section we derived a analytical tool to characterise the persistence of a population facing shifting range in a stochastic environment. We also analysed further the dependence of the criterion with respect to the parameters of the model and investigated the existence of a critical shifting speed $c$, separating persistence from extinction, in the particular case of Gaussian dispersal kernel.

\section{Numerical calculation of persistence conditions, with application to shifting butterfly populations}\label{sec:example}
We now focus on applying the stochastic model and theory to butterfly populations responding to climatic change. Our goal is twofold, first to illustrate the calculation of $\Lambda$ \eqref{Lambdar0K} and the critical shifting speed $c^*$ \eqref{eq:cstar} in a ecologically realistic problem and second to investigate the impact of variability in environmental shift ($\sigma$) and growth rate ($r$) on population persistence. 

We assume that for each generation, there are only two possible environments, a good environment when $\alpha=\alpha_g$ and a bad environment when $\alpha=\alpha_b$, such that 
$$P(\alpha=\alpha_b)=0.5 \text{ and } P(\alpha=\alpha_g)=0.5.$$ 
We consider $\sigma\in\{\underline{\sigma},\overline{\sigma}\}$, and $r\in\{\underline{r},\overline{r}\}$ with
\be\underline{\sigma}<0<\overline{\sigma}\text{ and } 0<\underline{r}<\overline{r}.\ee
To define good and bad environments we assume that the larger $r$, the better the environment and the smaller $\sigma$ the better the environment and thus consider
\be\alpha_g=(\underline{\sigma},\overline{r}) \text{ and } \alpha_b=(\overline{\sigma},\underline{r}).\ee
\citet{Lerouxetal13} study the effect of climate change and range migration for twelve species of butterflies in Canada. In their paper the authors use a reaction-diffusion model to study the ecological dynamics of the population and its persistence properties. In their analysis they estimate the shifting speed $c$ and the growth rate $r$ for each species to compute the critical diffusion coefficient that will determine the persistence of the population assuming no restriction on the length of the patch. Here we use their estimated means and standard deviation of $c$ and $r$ to fix the values of $c$, $\underline{\sigma}$, $\overline{\sigma}$, $\underline{r}$ and $\overline{r}$ defined above and
$$c=3.25\text{ km/year},\:\underline{\sigma}=-1.36\text{ km/year},\: \overline{\sigma}=1.36\text{ km/year},\:\underline{r}=2.07,\:\overline{r}=4.85$$
unless otherwise stated. 
We then study the persistence of the population in our stochastic framework for a fixed patch size $|\Omega_0|=10$ km. 

\subsection{Critical shifting speed for Gaussian kernel}\label{sec:numsimpower}
First, we investigate numerically the variation of the critical shifting speed $c^*$ as a function of the variance of the dispersal kernel in the case of Gaussian dispersal kernel. Indeed, using the analysis of Section \ref{sec:approxlambda}, we know that for Gaussian kernel, when $\sigma\equiv0$, that is the shifting speed is not stochastic anymore, then there exists a critical speed $c^*\geq0$ such that for all $c<c^*$ the population persists in the moving frame whereas if $c>c^*$ the population goes extinct. Using equation \eqref{eq:cstar} we know that $c^*=\sqrt{2(\sigma^K)^2\ln(\lambda_0\cdot \overline{R})},$ with $\lambda_0$ being the principal eigenvalue of the linear operation $\mathcal{K}_0$ defined in \eqref{Kc} when $c=0$ and $\overline{R}$ being the geometric average of the geometric growth rate, defined in \eqref{Rbar}. Notice from the derivation of $c^*$ in the previous section that when 
$\lambda_0\cdot\overline{R}<1$, the critical speed is zero because for all $c\geq0$, $\Lambda=e^{-\frac{c^2}{2(\sigma^K)^2}}\lambda_0\cdot \overline{R}$ defined in \eqref{LambdaGauss}, is always smaller than 1.
\begin{figure}
\begin{center}
\subfloat[]{\includegraphics[width=6.5cm, height=5cm]{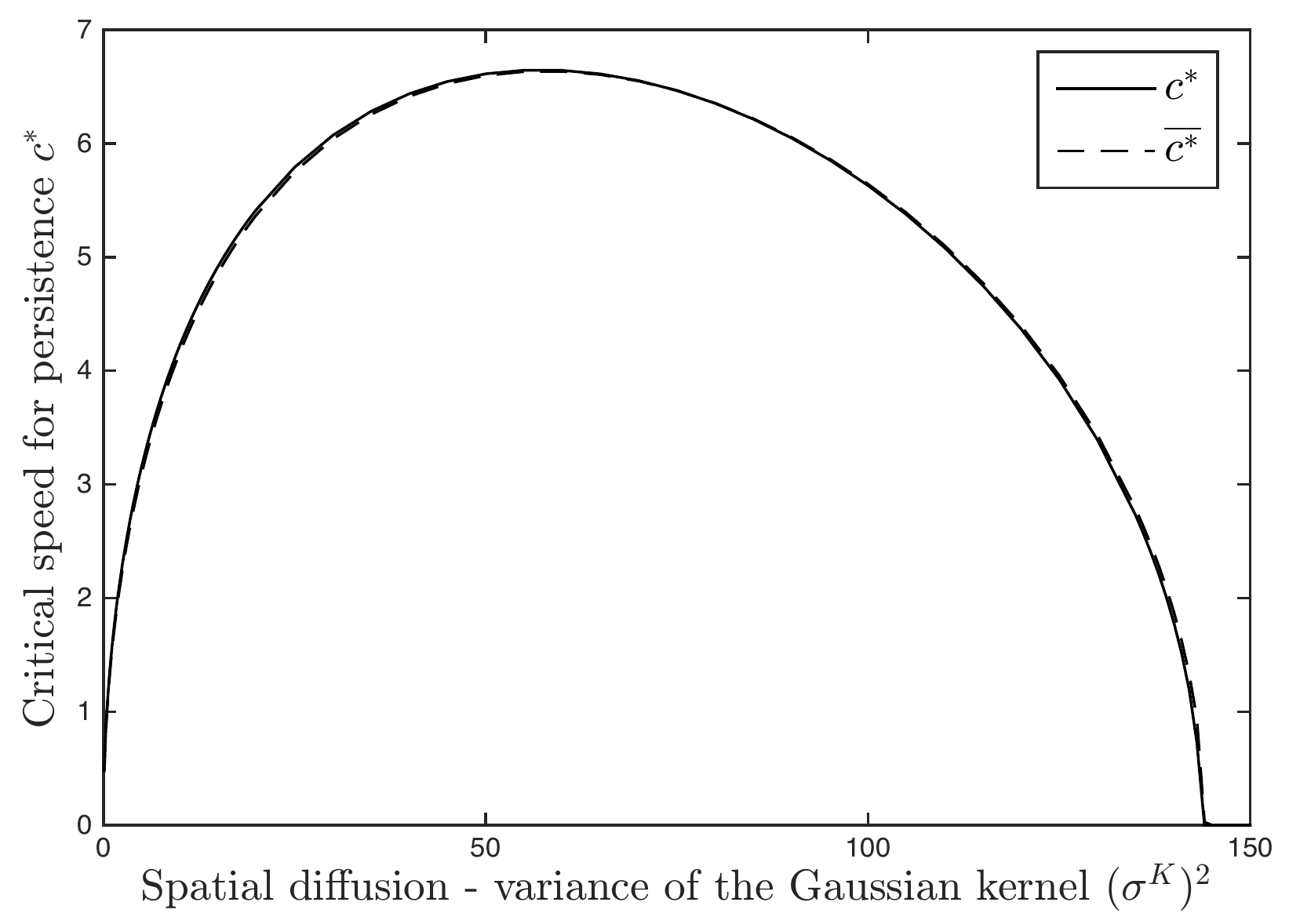}}\quad\quad
\subfloat[]{\includegraphics[width=6.5cm, height=5cm]{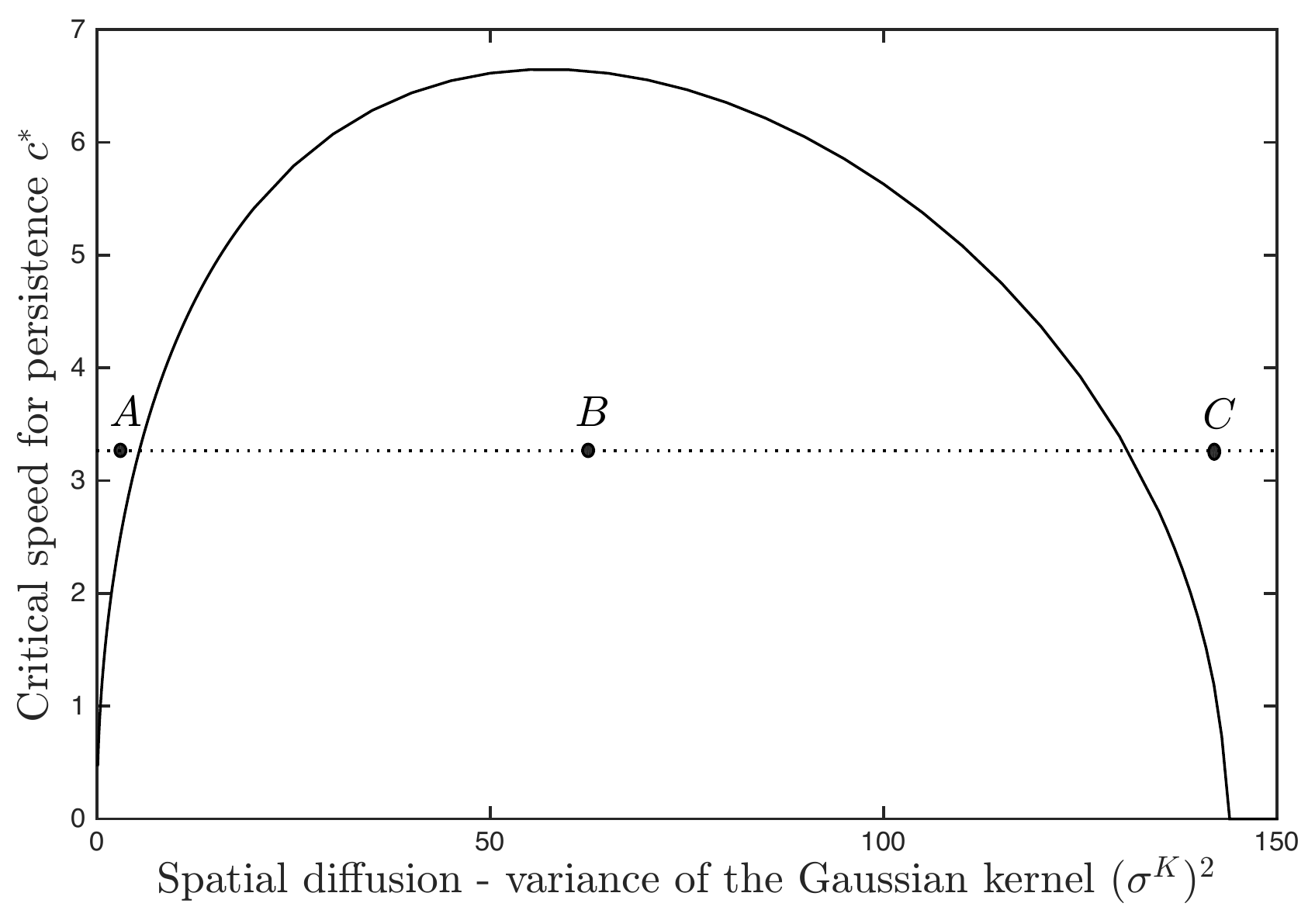}}
\caption{Critical speed for persistence $c^*$ as defined in \eqref{eq:cstar} and its approximation $\overline{c^*}$ using the dispersal success approximation $\overline{\lambda_0}$ instead of the principal eigenvalue $\lambda_0$, as a function of the dispersal Gaussian kernel variance (a). On the left panel (b), the critical speed for persistence $c^*$ is compared to a fixed shifting speed $c=3.25$ km/year. Points $A$, $B$ and $C$ highlight the existence of three different persistence regime as the variance of the dispersal kernel increases. The parameter values are the following:
$\underline{\sigma}= \overline{\sigma}=0\text{ km/year},\:\underline{r}=2.07,\:\overline{r}=4.85.$ and $|\Omega_0|=10$ km}\label{fig:CritSpeedGauss}
\end{center}
\end{figure}
As can be observed from Figure \ref{fig:CritSpeedGauss}, the critical shifting speed for persistence first increases with the variance of the dispersal kernel as the population is more and more mobile and can track its favorable environment more easily. Then the critical speed decreases for large values of $(\sigma^K)^2$ and converges to 0 as the patch size becomes too small for the population to persist even in a non-shifted environment. When the shifting speed $c=3.25$ km/year, as it was suggested by \citet{Lerouxetal13}, there exists three different regimes for the persistence of the population depending on the value of the variance $(\sigma^K)^2$ (Figure \ref{fig:CritSpeedGauss}b). When the variance of the dispersal kernel is small (point $A$), the critical speed for persistence $c^*$ is smaller than the shifting speed $c=3.25$ and thus population can not keep pace with its environment. As the variance $(\sigma^K)^2$ increases (point $B$), the critical speed for persistence increases above $3.25$ and thus the population persists. Nevertheless, when $(\sigma^K)^2$ becomes too large (compared to the patch size), the critical speed for persistence decreases below the shifting speed $c=3.25$ to reach 0 (point $C$) and the population does not keep pace with its environment anymore. One can also note that in this framework the dispersal success approximation for the principal eigenvalue $\overline{\lambda_0}$, defined in \eqref{eq:DispSuccAppr}, gives an accurate approximation of the critical speed for persistence (Figure \ref{fig:CritSpeedGauss}a).

\subsection{Effect of the dispersal on persistence of butterfly in Canada}
\begin{figure}
\begin{center}
\includegraphics[scale=0.5]{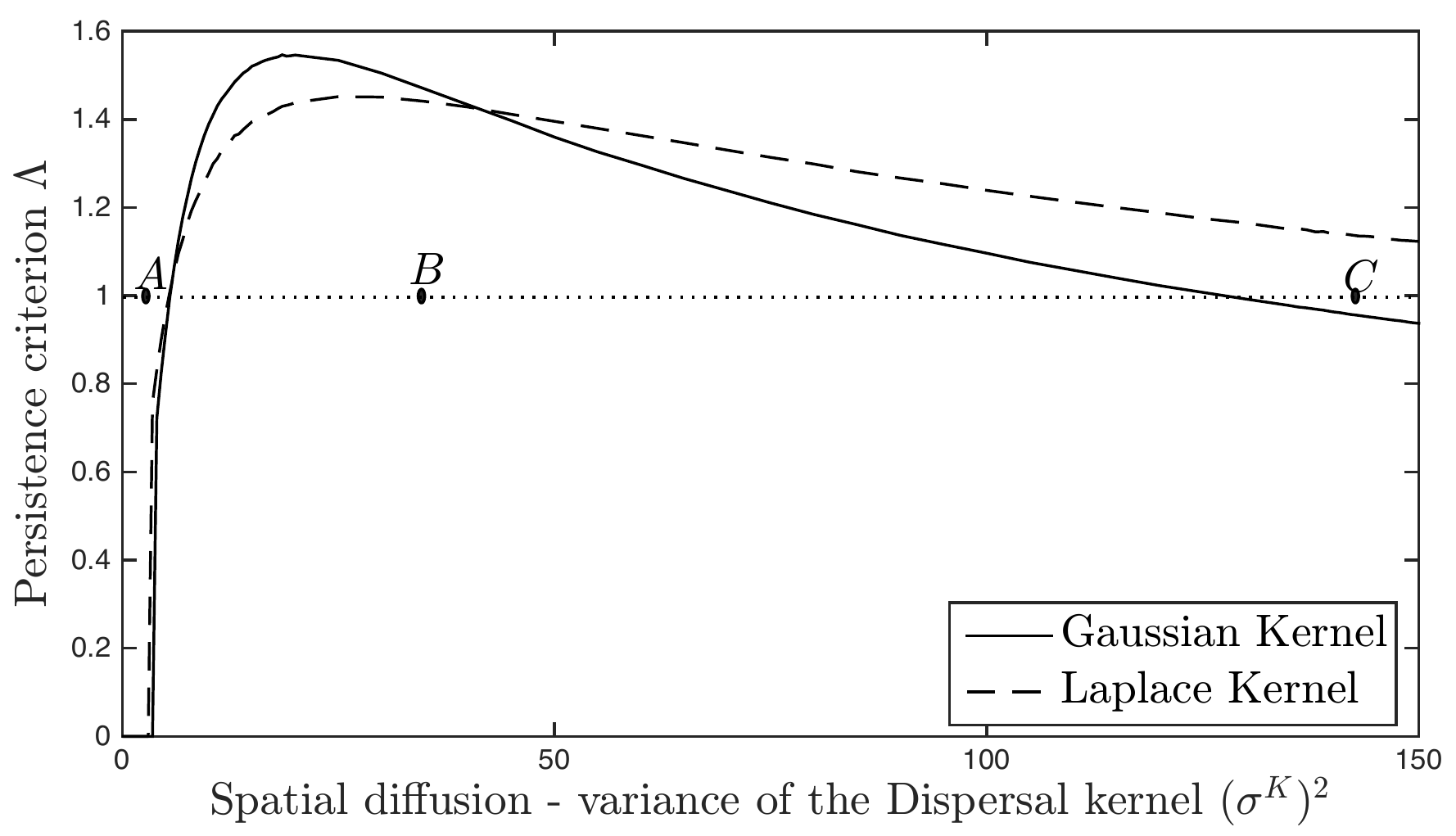}

\caption{Approximation of $\Lambda$ as a function of the dispersal kernel variance, for Gaussian kernel (solid curve) and Laplace kernel (dashed curve). The value of $\Lambda$ is compared to 1 (dotted curve). Points $A$, $B$ and $C$ highlight the existence of three different persistence regimes as the variance of the dispersal kernel increases. The parameter values are the following:
$c=3.25\text{ km/year},\: \underline{\sigma}=-1.36\text{ km/year},\: \overline{\sigma}=1.36\text{ km/year},\:\underline{r}=2.07,\:\overline{r}=4.85$ and $|\Omega_0|=10$ km }\label{fig:Lambda}
\end{center}
\end{figure}

In this section, we approximate the persistence criteria $\Lambda$, defined in \eqref{Lambda}, computing $\Lambda_t$ for large $t$ and study the change in the persistence criterion as a function of the variance of the dispersal kernel, when all the other parameters are fixed. We choose parameters based on the analysis by \citet{Lerouxetal13}: 
$$c=3.25\text{ km/year},\: \underline{\sigma}=-1.36\text{ km/year},\: \overline{\sigma}=1.36\text{ km/year},\:\underline{r}=2.07,\:\overline{r}=4.85.$$
We use two different kernels, the Laplace kernel
\be K(x)=\frac{\alpha}{2}e^{-\alpha|x|}\ee
and the Gaussian kernel
\be K(x)=\frac{1}{\sqrt{2\pi(\sigma^K)^2}}e^{-\frac{x^2}{2(\sigma^K)^2}}\ee
and consider for each kernel a variance $(\sigma^K)^2$, ($(\sigma^K)^2=(2/\alpha^2)$ for Laplace kernel) that spans from 0.1 to 150 $\text{km}^2$/year.
We approximate $\Lambda$, computing 
\be \Lambda_t=\left(\int_\Omega \tilde{n}_t(x)dx\right)^{1/t}\ee
for large $t$, as defined in \eqref{Lambda} (Figure \ref{fig:Lambda}). 
As can be expected from the computation of the critical speed in Figure \ref{fig:CritSpeedGauss}, the value of $\Lambda$ increases with the variance of the dispersal kernel at first and then decreases when the variance becomes too large with respect to the patch size (Figure \ref{fig:Lambda}). Whereas the persistence criteria $\Lambda$ associated with the Laplace kernel stay above 1 for large value of $(\sigma^K)^2$, the one associated with the Gaussian kernel decreases below 1 for large variance (Point $C$ in Figure \ref{fig:Lambda}). 

\subsection{Stochasticity of the parameters and persistence of the population}
We now turn our attention to investigating the effect of the stochasticity in the parameters on population persistence. We analyse the effects of increasing the variance of the yearly shift $\sigma$ or increasing the variance of the growth rate $r$ separately. In each case we assume that the expectation of the random variable $\sigma$ or $r$ are fixed and the probability of a bad (respectively good) environment is also fixed to 0.5. In this section we used the values from \citet{Lerouxetal13} to fix the parameters values and assume that 
$$c=3.25\text{ km/year}, \:|\Omega_0|=10\text{ km} \text{ and }(\sigma^K)^2=25 \text{ km}^2\text{/year}.$$

\begin{figure}
\begin{center}
\subfloat[]{\includegraphics[scale=0.4]{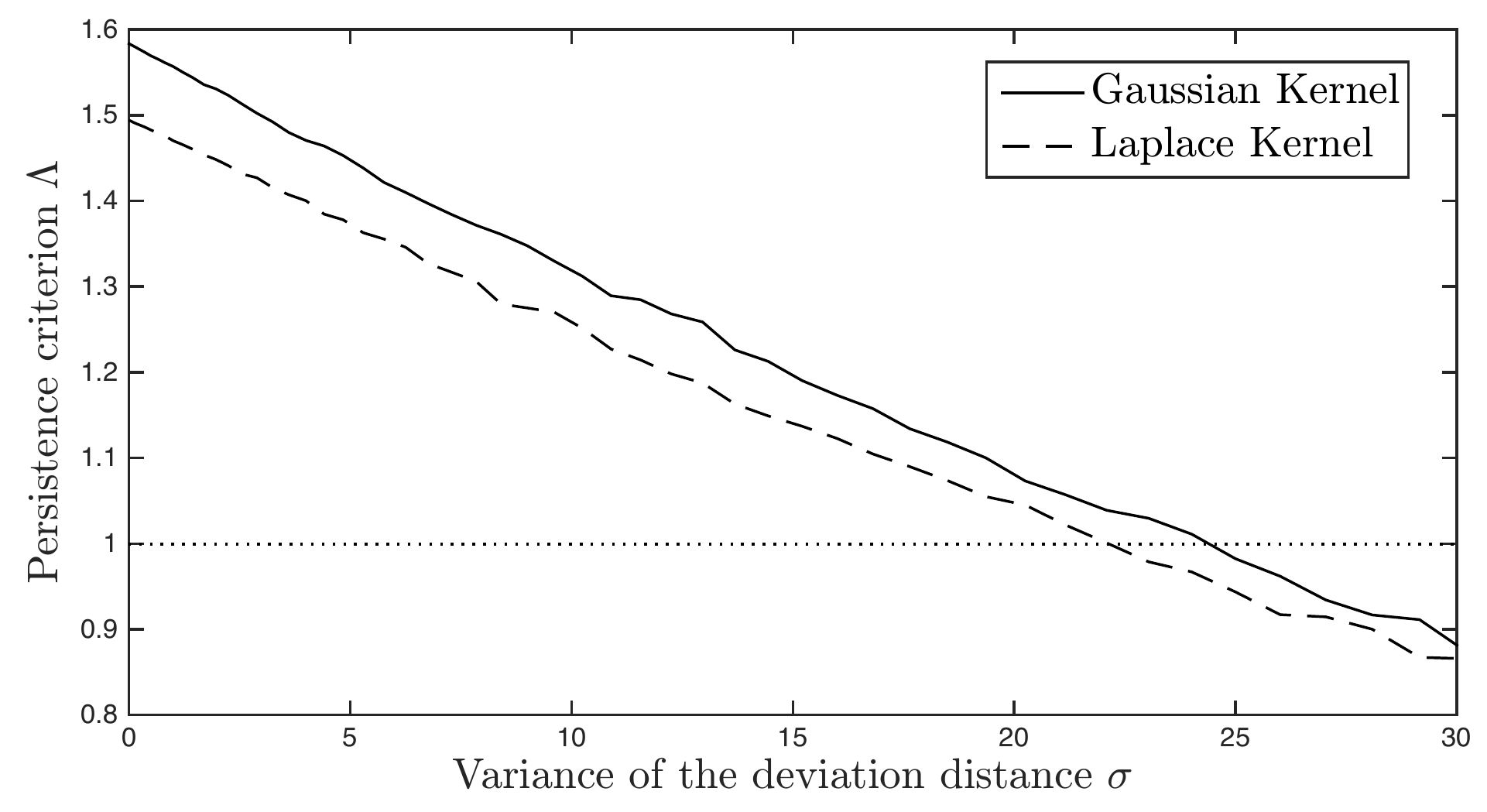}}\quad
\subfloat[]{\includegraphics[scale=0.4]{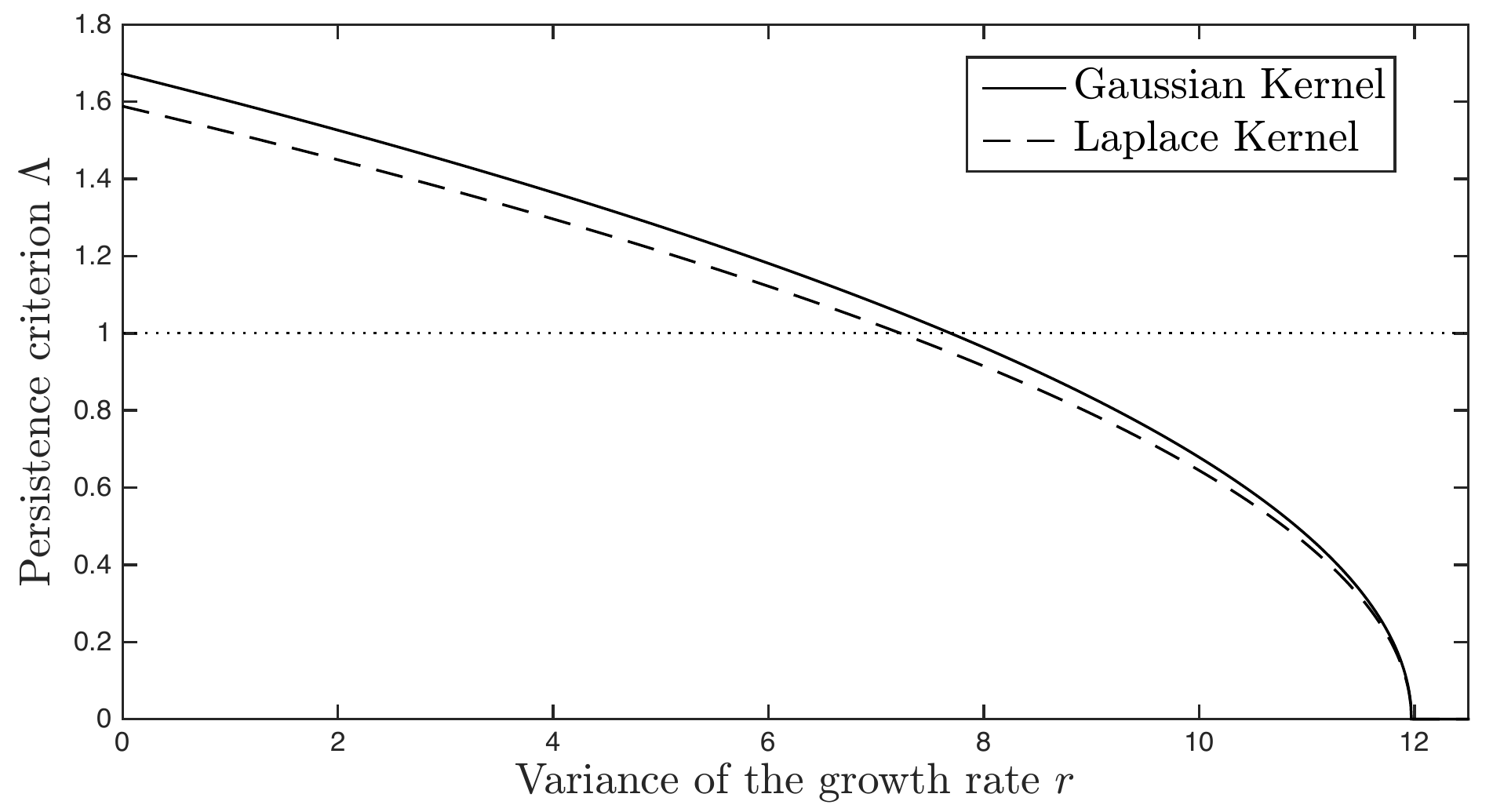}}

\caption{Approximation of $\Lambda$ as a function of the variance of the random shift $\sigma$ (panel (a)) or as a function of the variance of the growth rate $r$ for Gaussian Kernel (solid line) and Laplace kernel (dashed line). The parameter value for these computations are $c=3.25\text{ km/year}$, $|\Omega_0|=10$ km and the variance of the dispersal kernel is $25 \text{ km}^2\text{/year}$. For the left panel we fixed $\underline{r}=2.07,\:\overline{r}=4.85$, and for the right panel we fixed $\underline{\sigma}=-1.36,\:\overline{\sigma}=1.36$.}\label{fig:VarSigmaGrowth}
\end{center}
\end{figure}

We first analyse the variation of the persistence criterion $\Lambda$ as a function of the variance of the yearly shift $\sigma$. To do so we fix $\underline{r}$, $\overline{r}$ using again the values from \citet{Lerouxetal13} and assume that 
$$\underline{r}=2.07 \text{ and } \overline{r}=4.85.$$
We approximate $\Lambda$, computing $\Lambda_t$ for large $t$, as a function of the variance of $\sigma$, the yearly shift, assuming that $E[\sigma]=0$ (Figure \ref{fig:VarSigmaGrowth}(a)). For this analysis we consider again two types of dispersal kernel: a Gaussian dispersal kernel and a Laplace dispersal kernel. As illustrated in Figure \ref{fig:VarSigmaGrowth}(a), the persistence criterion $\Lambda$ decreases as the variance of $\sigma$ increases, and thus persistence is harder to achieve in environment with high stochasticity. Moreover when the variance of $\sigma$ becomes too large, the population can not persist anymore ($\Lambda<1$ in Figure \ref{fig:VarSigmaGrowth}(a)).

Using a similar approach we analyse the effect of increasing the variance of the growth rate $r$ on the persistence criterion $\Lambda$. To do so we fix $\underline{\sigma}$ and $\overline{\sigma}$ using again the values from \citet{Lerouxetal13} and assume that 
$$\underline{\sigma}=-1.36 \text{ and } \overline{\sigma}=1.36.$$
As illustrated in Figure \ref{fig:VarSigmaGrowth}(b), the persistence criterion is also decreasing as a function of the variance of the growth rate $r$. We thus conclude that adding variance to a model will impact negatively the persistence of the population.

In this section we considered a simple example, applicable to butterfly populations, where there exists only two outcomes for the random variable: a bad or a good environment. We investigated the dependence of the persistence of the population on the variance of the dispersal kernel and concluded that at first the persistence criterion (and also the critical shifting speed in the particular case of Gaussian kernels) increases with the variance of the dispersal kernel as the population increases its mobility and thus its ability to follow its favourable habitat. On the other hand the persistence criterion decreases as the variance of the dispersal kernel becomes too large, because the population flees its favourable habitat. We also investigated the dependence of the persistence of the population on the variability of the two stochastic variables in our model. We numerically analysed the persistence criterion as a function of the variance of the random variable $\sigma$, the yearly shift, or as a function of the variance of the random variable $r$, the growth rate. We observed that the persistence decreases with the variability of $\sigma$ (Figure \ref{fig:VarSigmaGrowth}(a)) and with the variability of $r$ (Figure \ref{fig:VarSigmaGrowth}(b)) and concluded that the variability of the asymptotic shifting speed or of the growth rate has a negative effect on population persistence.

\section{Discussion}
In this paper we model and analyse persistence conditions for a population whose range distribution is shifted toward the pole at a stochastic speed. The problem arises from the effects of climate change on range distribution \citep{P06}. While previous mathematical analyses of this problem consider deterministic environment \citep{ZK1,ZK2,HZHK, PK15}, it has been highlighted that the presence of climate change increases the temporal variability of the environment, partly due to extreme climatic events \citep{IPCC07}. Therefore, understanding the effect of a random geometric growth rate and a random shifting speed of the favorable environment becomes necessary. 

We first summarize the main results of our study. 
\begin{itemize}
\item Convergence in distribution to an equilibrium (as time goes to infinity), when the population is assumed to have compensatory competition dynamics,
\item Existence of a deterministic metric, having a biological interpretation in terms of a growth rate, which characterises persistence of the population in the case of compensatory and overcompensatory competition dynamics, 
\item Derivation of the persistence metric using the principal eigenvalue of a linear integral operator when the deviation of the shifting speed is fixed, 
\item Existence of a critical shifting speed separating persistence from extinction in the specific case of Gaussian dispersal kernel, 
\item Approximation of the critical shifting speed, using the dispersal success approximation \citep{VL97} and the modified dispersal success approximation \citep{RBM15}, 
\item Analysis of the persistence of a population of butterfly facing climate change and range shift as a function of its dispersal capacity, through numerical simulation,
\item Analysis of the qualitative effect of the variability in the yearly shift or in the growth rate on the persistence of a population in the particular case of butterfly.
\end{itemize} 

More precisely, our model assumes that the geometric growth rate for low density in the favorable environment is random from one generation to the next and that the bounded favorable environment is shifted at some asymptotic speed $c$, with an additional random yearly shift $\sigma_t$. The center of the favourable environment, $s_t$, is thus given by $s_t=ct+\sigma_t$. We use the theory of integrodifference equations to analyse the persistence of the population in this framework and employ a change of variables to study the problem in the shifted framework to track the favorable environment.  

The theory of \citet{Hardinetal1} and \citet{JJL} applies to our shifted problem and provides a general measure for persistence in term of the linearised problem at zero. Indeed in these two papers the authors prove that the persistence of the population depends on the magnitude of the metric $\Lambda$ \citep{JJL} defined as the asymptotic growth rate of the linearised operator, or equivalently on the magnitude of the metric $R$ \citep{Hardinetal1} defined as the asymptotic $L^\infty$ norm of the linearised operator. The former metric defined by \citet{JJL} has more biological meaning and is easier to compute numerically, which is why we chose to use it in this paper. Our paper thus provides a mathematical metric to measure persistence of the population facing climate change and shifted range distribution in a temporally variable environment. Note that we chose to study the problem in one dimension (i.e $\Omega\subset\R$) but the theory can be extended to $\Omega\subset\R^n$, $n\geq 1$, and the results in section \ref{sec:pers} are straightforward to extend.

We use the analytical results of Section \ref{sec:pers} to study the persistence of a butterfly population facing range shifts. We computed, for two different dispersal kernel (Laplace and Gaussian), the persistence metric as a function of the dispersal capacity (variance of the dispersal kernel). We find that for a fixed shifting speed (estimated from a paper by \citet{Lerouxetal13}), there exists three different regimes (Figure \ref{fig:Lambda}). When the dispersal capacity of the population is small, the butterflies can not keep track with their favorable environment and go extinct, when the dispersal capacity increases the persistence metric increases above one and the population persists. On the other hand, as the dispersal variance becomes too high, the population disperses outside its favourable habitat and goes extinct.

In the special case where the shift is deterministic, that is the yearly variation is $\sigma_t=0$ for all $t\in\N$, then the problem describes the dynamics of a population with variable growth, facing shifting ranges. In this case we characterise the persistence of the population through the magnitude of the asymptotic growth rate of a linear operator, similarly to the results by \citet{ZK1} in the case of deterministic scalar integrodifference equations with shifting boundaries. We also show that in the case of Gaussian dispersal kernel, this asymptotic growth rate is decreasing with respect to $c$. This proves the existence of a critical speed $c^*$ such that for all shifting speed $c$ less than $c^*$, the asymptotic growth rate will be above one and the population will persist whereas when the shifting speed is above $c^*$, the asymptotic growth rate is below one and the population converges to zero in its favorable environment. These results are also independent of the initial condition, as soon as it is non negative and non trivial and in this sense is similar to the one proved by \citet{BDNZ} for deterministic scalar reaction-diffusion equations.

To analytically determine the critical speed for persistence, $c^*$, we need to compute the principal eigenvalue of the linearised operator around the trivial steady state. As the dispersal kernel is positive everywhere, this principal eigenvalue can be computed using the power method. Nevertheless in a biological framework, this principal eigenvalue can be approximated using the dispersal success function \citep{VL97} or a modified dispersal success function \citep{RBM15} (Figure \ref{ppleeigenDispSuccAppr}). In this case the principal eigenfunction is either approximated by the probability to stay in the suitable habitat, after the dispersal stage, normalised by the size of the habitat (dispersal success approximation), or by the probability to stay within the favourable habitat, weighted by the proportion of individual at each point of the favourable domain, normalised again by the size of the habitat \citep{RBM15}. As illustrated in Figure \ref{fig:CritSpeedGauss}a, the dispersal success approximation gives accurate results for the computation of the critical speed for persistence. On the other hand, these approximations do not give accurate estimates of the principal eigenfunction when the dispersal kernel is shifted by $c$ (if one wanted to estimate $\lambda_c$ directly instead indirectly through $\lambda_0$). In these case one could use the methods described by \citet{PK15} to compute efficiently the principal eigenvalue of a given linear operator.

We also highlighted, using numerical simulation, the negative effect of the variability of the parameter on the persistence of the population. Indeed, as illustrated in Figure \ref{fig:VarSigmaGrowth}, the persistence metric decreases with the variance of both parameters: $\sigma$, the yearly shift and $r$, the growth rate. 

Finally, note that the asymptotic shifting speed $c$ is actually uncertain in the sense that it depends on the severity of the climate change \citep{IPCC14}. In this paper we assumed that $c$ was known but if one would only know the distribution of the different outcomes for $c$, the asymptotic shifting speed, one can then compute the probability of persistence using our analysis. Indeed, in the specific case of Gaussian kernel and assuming that the yearly shift $\sigma_t$ is zero for all $t$, then the probability of persistence would be given by the probability that $c$ is less than $c^*$. 

\section*{Aknowledgments}
MAL gratefully acknowledges funding from aCanada Research Chair and on NSERC Discovery Grant. JB gratefully acknowledges funding  for a postdoctoral fellowship from the Pacific Institute for Mathematical Sciences (PIMS).

\appendix
\section*{Appendix}
\section{Proof of Theorem \ref{th:convergence} and \ref{th:persistence}}\label{sec:proofhardinetal}
%\subsection{Proof of Theorem \ref{th:convergence} and \ref{th:persistence}}\label{sec:proofhardinetal}
The proof of Theorems \ref{th:convergence} and \ref{th:persistence} mostly follows from \citet{Hardinetal1}[Theorem 4.2 and Theorem 5.3]. Indeed \citet{Hardinetal1} study the general stochastic model
\be X_{t+1}=F_{\alpha_t}(X_t),\ee
with $(X_t)_t$ a random process that takes values in the set on non negative continuous function in $\Omega$, $(\alpha_t)_t$, independent identically distributed random variable taking values in the set $\mathcal{S}$. And they have the following assumption for $F_\alpha$
\begin{itemize}
\item[(H1)]  For each $\alpha\in\mathcal{S}$, $F_\alpha$ is a continuous map of $C_+(\Omega)$, into itself such that $F_\alpha(u)=0\in C_+(\Omega)$ if and only if $u=0\in C_+(\Omega)$.
\item[(H2)] If $u,\: v\in C_+(\Omega)$ and $u\geq v$ then $F_\alpha(u)\geq F_\alpha(v)$.
\item[(H3)] There exists some $b>0$ such that for $u\in C_+(\Omega)$
\begin{itemize}
\item[(a)] $||F_\alpha(u)||_\infty\leq b$ for all $\alpha\in\mathcal{S}$ and $||u||_\infty\leq b$,
\item[(b)] there exists some time $t$ (depending on $u_0$) such that 
\be ||F_{\alpha_t}\circ \cdots \circ F_{\alpha_0}u_0||_\infty<b,\ee
for all $\alpha_0,\dots,\alpha_t\in\mathcal{S}$,
\item[(c)] there exists some $d>0$ such that 
\be F_\alpha(b)\geq d\ee
for all $\alpha\in\mathcal{S}$.
\end{itemize}
\item[(H4)] Let $B_b:=\{u\in C_+(\Omega):\: ||u||_\infty\leq b\}$, then there is some compact set $D\subset C_+(\Omega)$ such that $F_\alpha(B_b)\subset D$ for all $\alpha\in\mathcal{S}$.
\item[(H5)] There exists some $h>0$ such that 
\be ||F_\alpha(u)||_\infty\leq h||u||_\infty,\ee
for all $\alpha\in\mathcal{S}$ and $u\in C_+(\Omega)$.
\item[(H6)] There exists some $\xi>0$ such that 
\be F_\alpha(B_b)\subset K_\xi,\ee
where $K_\xi:=\{u\in C_+(\Omega):\:u\geq \xi||u||_\infty\}$.
\item[(H7)] For each $a>0$, there exists a continuous function $\tau:(0,1]\to(0,1]$ such that $\tau(s)>s$ for all $s\in(0,1)$ and such that 
\be \tau(s)F_\alpha(u)\leq F_\alpha(su)\ee
for all $\alpha\in\mathcal{S}$ and $u\in C_+(\Omega)$ such that $a\leq u\leq b$.
\item[(H8)] $F_\alpha$ is Fr\'echet differentiable (with respect to $C_+(\Omega)$) at $0\in C_+(\Omega)$. We denote by $\mathcal{L}_\alpha$ the operator $F'_\alpha(0)$
\item[(H9)] There exists a function $\mathcal{N}:\R_+\to[0,1]$ such that 
\be \underset{u\to0,\:u>0}{\lim} \mathcal{N}(u)=1 \text{ and } \mathcal{N}(||u||_\infty)\mathcal{L}_\alpha u\leq F_\alpha(u)\leq \mathcal{L}_\alpha u,\ee
for all $u\in C_+(\Omega)$.
\end{itemize}
We want to prove that the previous hypotheses are satisfied in our framework and then apply Theorem 4.2, Theorem 5.3 by \citet{Hardinetal1} and Theorem 2 by \citet{JJL}.
One can check, as it is done by \citet[Section 5.3]{JJL}, that under Hypotheses \ref{hypK}-\ref{alphastar}, the previous hypotheses (H1)-(H9) are satisfied and one can then apply Theorem 4.2, Theorem 5.3 from \citet{Hardinetal1} and Theorem 2 from \citet{JJL}. For completeness we will write the main steps referring most of the time to \citet[Section 5.3]{JJL}. First, we denote by, 
\be\label{Kbar}\overline{K}:=\sup\{K(x), x\in\R\}<+\infty,\ee
and 
\be\label{Kubar}\underline{K}:=\inf\left\{ K(x), \: x\in( \inf \Omega-\sup\Omega+c,\sup \Omega-\inf\Omega+c)\right\}>0.\ee
These two constants will be used several times in the proof of (H1)-(H9) below.  

\begin{itemize}
\item[(H1)] The continuity of $F_\alpha$ follows from the continuity of $f_r$, and the boundedness of $K$ and $\int_\Omega g_0(y)dy$. $K$ positive in $\R$, $f_r$ and $g_0$ non negative yield the second statement.
\item[(H2)] This follows from the monotonicity of $f_r$ for all $\alpha\in\mathcal{S}$.
\item[(H3)] The constant $b>0$ will be defined later in the proof,
\begin{itemize}
\item[(a)] for all $\alpha\in\mathcal{S}$, $u\in C_+(\Omega)$,
\begin{align}
||F_\alpha(u)||_\infty&=\underset{x\in\Omega}{\max}\int_\Omega K(x-y+c)g_0(y-\sigma)f_r(u(y))dy\\
&< m \cdot \overline{K}\cdot\int_{\Omega_0} g_0(y)dy=b.
\end{align}
This proves the statement.
\item[(b)] Using the same argument as before, for all $u_0\in C_+(\Omega)$, $t\in\N^*$, $\alpha_t,\dots,\alpha_0\in\mathcal{S}$,
\be ||F_{\alpha_t}\circ \cdots\circ F_{\alpha_0}(u_0)||_\infty=||F_{\alpha_t}(u)||_\infty<b,\ee
where $u\in C_+(\Omega)$.
\item[(c)] For all $x\in\Omega$,
\be F_\alpha(b)(x)\geq \underset{\alpha\in\mathcal{S}}{\inf}\:f_r(b)\int_\Omega K(x-y+c)g_0(y-\sigma)dy\geq \underset{\alpha\in\mathcal{S}}{\inf}\:f_r(b)\cdot d_1\ee
with 
\be d_1:= \underline{K}\cdot\int_{\Omega_0} g_0(y)dy>0\ee
One concludes using the positivity of $K$ and Hypothesis \ref{hypf}(iii)b that there exists $d>0$ such that for all $x\in\Omega$,
\be F_\alpha(b)(x)\geq d.\ee
\end{itemize}
\item[(H4)] This statement follows from the continuity of $K$, uniform boundedness of $F_\alpha$ and $f_r$ and Hypothesis \ref{alphastar}, details can be found in \citet{JJL}[Section 5.3]. 
\item[(H5)] From assumption $\ref{hypf}(ii)c$, we have that for all $u>0$, $f'_r(0)u\geq f_r(u)$. Thus, using this inequality , with assumptions \ref{hypK}(ii), \ref{hyp:Q}(ii) and \ref{hypf}(iii)a we get for all $u\in C_+(\Omega)$, for all $\alpha\in\mathcal{S}$,
\be ||F_\alpha(u)||_\infty\leq h||u||_\infty,\ee
with $h:=\overline{r}\cdot\overline{K}\cdot \int_\Omega g_0(y)dy$.
\item[(H6)] First, notice that for all $\alpha\in\mathcal{S}$, $u\in B_b$
\be ||F_\alpha(u)||_\infty\leq \overline{r}\cdot \overline{K}\cdot\int_\Omega g_0(y-\sigma)u(y)dy \: \implies\: \int_\Omega g_0(y-\sigma)u(y)\geq \frac{||F_\alpha(u)||_\infty}{\overline{r}\cdot\overline{ K}},\ee
and using Hypothesis \ref{hypf}(ii)c,
\be f_\alpha(u)\geq \frac{f_\alpha(b)}{b}u.\ee
Then for all $x\in\Omega$
\begin{align}
F_\alpha(u)(x)&\geq \underline{K}\cdot\frac{f_\alpha(b)}{b}\int_\Omega g_0(y-\sigma)u(y)dy\\
&\geq \frac{\underline{K}}{\overline{K}}\cdot\frac{f_\alpha(b)}{b}\frac{||F_\alpha(u)||_\infty}{\overline{r}}
\end{align}
and the statement is proved.
\item[(H7)] This proof is also derived from \citet{JJL}. We want to find a continuous function $\tau:\: (0,1]\to(0,1]$ such that for all $s\in(0,1)$,
\begin{align*}
&\tau(s)F_\alpha(u)\leq F_\alpha(su)\\
\Leftrightarrow & \int_\Omega K(x-y+c)g_0(y-\sigma)\tau(s)f_r(u(y))dy\leq \int_\Omega K(x-y+c)g_0(y-\sigma)f_r(su(y))dy
\end{align*}
Thus it is sufficient to have $\tau(s)f_r(u(y))\leq f_r(su(y))$ for all $y\in\Omega$. From Hypothesis \ref{hypf}(ii)c, as $s\in(0,1)$, we have that for all $\alpha\in\mathcal{S}$
\be s<\frac{f_r(su)}{f_r(u)},\ee
and letting $\tau(s)=\min\left\{\frac{f_r(su)}{f_r(u)},\: \alpha\in\mathcal{S},\: a\leq u\leq b\right\}$ we have, for all $s\in(0,1)$, $\alpha\in\mathcal{S}$ and $u\in[a,b]$,
\be \tau(s)>s \text{ and } \frac{f_r(su)}{f_r(u)}\geq \tau(s)\ee
and the statement is proved.
\item[(H8)] We want to prove that 
\be \underset{h\to0}{\lim}\frac{||F_\alpha(0+h)-F_\alpha(0)-\mathcal{L}_\alpha h||_\infty}{||h||_\infty}=0,\ee
where $\mathcal{L}_\alpha:=F'_\alpha(0)$ is a linear operator. Using the differentiability of $f_r$ at 0, one proves that the limit exists and $\mathcal{L}_\alpha h=r\int_\Omega K(x-y+c)g_0(y-\sigma)h(y)dy$.
\item[(H9)] The second part of the inequality follows from assumption \ref{hypf}(ii)c. For each $\alpha\in\mathcal{S}$ let $\mathcal{N}_\alpha:[0,+\infty)\to[0,1]$ be such that 
\begin{equation}
\mathcal{N}(u)=\begin{cases} \frac{f_\alpha(u)}{r\cdot u} &\text{if } u>0,\\
1 &\text{if } u=0.\end{cases}\end{equation}
The function $\mathcal{N}_\alpha$ is continuous and defining $\mathcal{N}(u)=\min\left\{\mathcal{N}_\alpha(u),\: \alpha\in\mathcal{S}\right\}$, it gives the wanted statement (using the fact that $\mathcal{N}_\alpha(u(x))\geq \mathcal{N}_\alpha(||u||_\infty)$).
\end{itemize}
We can thus apply the Theorem 4.2 from \citet{Hardinetal1} to prove Theorem \ref{th:convergence} and use Theorem 5.3 by \citet{Hardinetal1} and Theorem 2 by \citet{JJL} to prove Theorem \ref{th:persistence}. For weak of completeness we include the three Theorems from \citet{Hardinetal1} and \cite{JJL} below.
\begin{theorem} \textbf{ - [\citet[Theorem 4.2]{Hardinetal1}]}\\
Suppose that Hypotheses (H1)-(H7) above and  Hypothesis \ref{alpha}(i) in Section \ref{sec:model} are satisfied and that $u_0\neq0\in C_+(\Omega)$ with probability one. Then $u_t$ solution of \eqref{IDEmoving} converges in distribution to a stationary distribution $\mu^*$, independent of $u_0$, such that either $\mu^*(\{0\})= 0$ or $\mu^*(\{0\}) = 1$.
\end{theorem}
\begin{theorem} \textbf{ - [\citet[Theorem 5.3]{Hardinetal1}]}\\
Suppose that Hypotheses (H1)-(H9) above and  Hypothesis \ref{alpha}(i) in Section \ref{sec:model} are satisfied and that $u_0\neq0\in C_+(\Omega)$ with probability one. Let $\mu^*$ be as in Theorem 4.2 and define $R=\lim_{t\to+\infty} ||\mathcal{L}_{\alpha_t}\circ\dots\circ\mathcal{L}_{\alpha_0}||$, with $\mathcal{L}$ the linearised operator around zero (\eqref{eq:utilde})
\begin{enumerate}[(a)]
\item If $R< 1$ then $\mu^*(\{0\}) = 1$ and $u_t\to0$ with probability one.
\item If $R > 1$ then $\mu^*(\{0\}) =0$.
\end{enumerate}\end{theorem}

\begin{theorem} \textbf{ - [\citet[Theorem 2]{JJL}]}\\
Let $R$ be defined as in the previous theorem by $R=\lim_{t\to+\infty}||\mathcal{L}_{\alpha_t}\circ\dots\circ\mathcal{L}_{\alpha_0}||$, then $\Lambda=R$ ($\Lambda$ defined in \eqref{Lambda}).\end{theorem}

\section{Persistence condition and spreading speed}\label{a:spreadspeed}
In this section we derive an heuristic criteria for persistence inspired from \citet{NKL00} linking the critical patch size and the asymptotic spreading speed to get necessary conditions for persistence. We consider problem \eqref{IDEnonmoving}, where the suitability functions $g_0$ are indicator functions, that is
\begin{equation*}g_0(y)=\mathds{1}_{\Omega_0}=\begin{cases} 1 &\text{if }y\in\Omega_0,\\
									0 &\text{otherwise.}\end{cases}\end{equation*}
									
In addition to the assumptions made in Section \ref{sec:model} we will also assume that $K$ is a thin tail dispersal kernel in the sense that it has exponentially bounded tails, i.e there exists $s>0$,
\begin{gather}\label{hyp:K2}
\int _\R e^{s|x|}K(x)dx<+\infty.
\end{gather}
This last assumption guarantees that the moment generating function (\eqref{MomGenFunc}) exists on some open interval of the form $(0,s^+)$. This assumption was not necessary to derive the persistence condition in Section \ref{sec:pers} but it will be used to study the speed of the stochastic wave. 
									
\subsection{Critical domain size in the constant environment}\label{a:stochspeed}
Let us first consider the problem in the non shifted frame and thus assume that $s_t\equiv0$, where $(s_t)_t$, defined in Section \ref{sec:model}, is the center of the suitable habitat at generation $t$. Thus we consider the following integrodifference equation
\be u^0_{t+1}(x)=\int_{\Omega_0} K(\xi-\eta)f_{r_t}(u^0_t(\eta))d\eta,\ee
denoting by $u^0$ the solution in this non shifted framework. Using the theorems in section \ref{sec:pers}, we get that as time goes to infinity, $u^0$ persists if 
\be\label{perscond0}E[\ln(r_0)] >-\ln(\lambda_0),\ee
where $E[\cdot]$ is the expectation of a random variable and $\lambda_0$ is the principal eigenvalue of the linear operator $\mathcal{K}_0$:
\be \mathcal{K}_0[u](x)=\int_{\Omega_0}K(x-y)u(y)dy.\ee
Now assume that $K$, $\Omega_0$ and $(r_t)_t$ are such that \eqref{perscond0} is satisfied and study the problem in the non moving, homogeneous framework.

\subsection{Spreading speed in a stochastic homogeneous environment}
Now we are interested in deriving the asymptotic spreading speed of the population in an homogeneous environment to compare it with the forced shifting speed $c$.
We thus consider the homogeneous problem on $\R$ in the non moving frame, i.e let $(n_t)_t$ be the solution of the equation
\be\label{homoIDE}
n_{t+1}(\xi)=\int_\R K(\xi-\eta)f_{r_t}(n_t(\eta))d\eta.
\ee
As we are considering the initial problem in the non moving frame, $K$ does not depend on $c$ and $g_0\equiv 1$ in $\R$ and thus the stochasticity comes only from the growth term. From the analysis of \citet{NKL00}, there are two different approaches to estimate the spreading speed of the stochastic process $(n_t)_t$. One can either consider the spreading speed of the expected wave or the asymptotic speed of the stochastic wave. We will only consider the latter approach and assume that the speed is governed by the linearisation at $0$. 
Denote by $(\tilde{n}_t)_t$, the solution of the linearised operator at 0, i.e for all $t\in\N$,
\be\label{linearhomoIDE}
\tilde{n}_{t+1}(\xi)=\int_\R K(\xi-\eta)r_t\tilde{n}_t(\eta)d\eta.\ee
We define the random variable $\Xi_t$ as the most rightward position such that $\tilde{n}_t$ is greater that some threshold, i.e
\be \Xi_t=\sup\left\{\xi\in\R,\: \tilde{n}_t>\overline{n}\right\},\ee
where $\overline{n}\in(0,1)$ is a fixed critical threshold. Assume that $\forall \xi\in\R$, $n_0(\xi)=\alpha e^{-s\xi}$, for some $s>0$, i.e the initial condition has a wave shape, then for all $t\in\N$, $\xi\in\R$,
\be \tilde{n}_{t+1}(\xi)=\alpha\prod_{i=0}^t(r_iM(s))e^{-s\xi},\ee
where $M$ is the moment generating function of $K$ (\eqref{MomGenFunc}). 
This function exists in some interval $(0,s^+)$ because of assumption \eqref{hyp:K2}. 
Moreover $\overline{n}=n_0(\Xi_0)=\tilde{n}_{t+1}(\Xi_{t+1})$, thus denoting by $\overline{c}_t(s)$ the spreading speed of $(\Xi_t)_t$ starting with $n_0(\xi)=\alpha e^{-s\xi}$ for all $\xi\in\R^+$, we have
\begin{align*}
\overline{c}_{t+1}(s)&=\frac{\Xi_{t+1}-\Xi_0}{t+1}\\
&=\frac{1}{t+1}\sum_{i=0}^t \frac{1}{s}\ln(r_iM(s))\\
&\left(=\frac{1}{s}\ln(M(s))+\frac{1}{t+1}\sum_{i=0}^t \frac{1}{s}\ln(r_i)\right).
\end{align*}
Thus $(\overline{c}_t(s))$ is the sum of independent identically distributed variables and thus converges in distribution to a random variable that is normally distributed with mean $\mu(s)$ and variance $\sigma^2(s)$ such that 
\be \mu(s)=E[\frac{1}{s}\ln(r_0M(s))]\ee
and \be \sigma^2(s)=\underset{t\to+\infty}{\lim}\frac{1}{t} V[\frac{1}{s}\ln(r_0M(s))]=0.\ee
As $\sigma^2(s)\equiv 0$, this implies that $\overline{c}_t(s)$ converges in probability to the constant $\frac{1}{s}E[\ln(r_0M(s))]$. This is true for all $s$ such that $M(s)$ exists. Now if we want to consider the more general cases when $n_0$ is a compactly supported function, the minimal speed over all the $s$ will be the relevant one and we have that the spreading speed of the stochastic wave at time $t$, $\overline{c}_t$, has mean $\mu=\underset{s>0}{\inf}\mu(s)$ and variance $(\sigma^*_t)^2=\sigma^2_t(s^*)$, where $s^*$ is such that $\mu(s^*)=\mu$, 
and thus converges in probability to 
\be \overline{c}^*=\underset{s>0}{\inf}\frac{1}{s}E[\ln(r_0M(s))].\ee

\bibliography{BLbiblio}
\bibliographystyle{spbasic}
\end{document}